\documentclass{article}

\usepackage{calrsfs}
\usepackage{amscd}
\usepackage{amsmath}
\usepackage{amsfonts}
\usepackage{amssymb}
\usepackage{amsthm}
\usepackage{hyperref}
\usepackage{comment}

\newtheorem{thm}{{\bf{\small T}{\scriptsize HEOREM}}}[section]
\newtheorem{prop}[thm]{{\bf{\small P}{\scriptsize ROPOSITION}}}
\newtheorem{lem}[thm]{{\bf{\small L}{\scriptsize EMMA}}}
\newtheorem{cor}[thm]{{\bf{\small C}{\scriptsize OROLLARY}}}
\newtheorem{defn}[thm]{{\bf{\small D}{\scriptsize EFINITION}}}
\newtheorem{rmq}[thm]{{\bf{\small R}{\scriptsize EMARK}}}
\newtheorem{example}[thm]{{\bf{\small E}{\scriptsize XAMPLE}}}

\newcommand{\MakeThmProofTitle}[1]{{\bf {\small{P}{\scriptsize ROOF
OF }{\small{T}}{\scriptsize HEOREM} #1}}}

\newcommand{\mum}{m}

\newcommand{\lawto}{\stackrel{\scriptscriptstyle{law}}{\longrightarrow}}

\renewcommand{\tilde}{\widetilde}

\newcommand{\norm}[1]{\left\| #1 \right\|}
    \DeclareMathOperator{\Jac}{Jac}

    \DeclareMathOperator{\Id}{Id}

    \DeclareMathOperator{\sgn}{sgn}

\newcommand{\dd}{\, {\rm d}}

\renewcommand{\geq}{\geqslant}
\renewcommand{\leq}{\leqslant}

\newcommand{\boE}{\mathcal{E}}
\newcommand{\E}{{\mathbb E}}

\newcommand{\N}{\mathbb{N}}
\renewcommand{\P}{\mathbb{P}}
\newcommand{\Z}{\mathbb{Z}}
\newcommand{\R}{\mathbb{R}}

\newcommand{\C}{\mathbb{C}}

\renewcommand{\phi}{\varphi}
\renewcommand{\epsilon}{\varepsilon}
\newcommand{\tq}{\, :\, }

\newcommand{\boL}{\mathcal{L}}

\newcommand{\boB}{\mathcal{B}}

\newcommand{\boN}{\mathcal{N}}
\newcommand{\boH}{\mathcal{H}}
\newcommand{\boF}{\mathcal{F}}

\newcommand{\boG}{\mathcal{G}}
\newcommand{\boW}{\mathcal{W}}

\newcommand{\RandVar}{Z}
\newcommand{\id}{\textup{i}}

\def\1{{\mathchoice {\rm 1\mskip-4mu l} {\rm 1\mskip-4mu l}
{\rm 1\mskip-4.5mu l} {\rm 1\mskip-5mu l}}}

\font\gfont=cmmi10 scaled \magstep{1.5}     
\newcommand{\gdelta}{\hbox{\gfont \char14}}

\begin{document}

\title{On almost-sure versions\\ of classical limit theorems\\
for dynamical systems}
\author{
J.-R. Chazottes\footnote{CPhT, CNRS-Ecole Polytechnique, 91128
Palaiseau Cedex, France, and CMM, UMI CNRS 2807, Universidad de
Chile, Av. Blanco Encalada 2120, Santiago, Chile,
jeanrene@cpht.polytechnique.fr}\\
S. Gou\"{e}zel\footnote{IRMAR, Universit\'{e} de Rennes 1, Campus de
Beaulieu, 35042 Rennes Cedex, France,
sebastien.gouezel@univ-rennes1.fr} }
\date{May 31, 2006}

\maketitle

\begin{abstract}
The purpose of this article is to support the idea that ``whenever
we can prove a limit theorem in the classical sense for a dynamical
system, we can prove a suitable almost-sure version based on an
empirical measure with log-average''. We follow three different
approaches: martingale methods, spectral methods and induction
arguments. Our results apply, among others, to Axiom A maps or flows, 
to systems inducing a Gibbs-Markov map, and to the stadium billiard.

\bigskip


\end{abstract}

\newpage

\tableofcontents

\newpage

\section{Introduction}

There has been recently a lively interest in probability theory
concerning almost-sure versions of classical limit theorems. The
prototype of such a theorem is the almost-sure central limit
theorem: if $\RandVar_n$ is an i.i.d.\ $L^2$ sequence with $\E(\RandVar_i)=0$ and
$\E(\RandVar_i^2)=1$, then, almost surely,
  \begin{equation}\label{ascltiid}
  \frac{1}{\log n} \sum_{k=1}^n \frac{1}{k} \gdelta_{\sum_{j=0}^{k-1}
  \RandVar_j/\sqrt{k}} \lawto \boN(0,1)
  \end{equation}
where ``$\lawto$'' means weak convergence of probability measures on
$\R$. Here and henceforth, $\delta_x$ is the Dirac mass at $x$. This
result should be compared to the classical central limit theorem,
which can be stated as follows:
\begin{equation}\label{cltiid}
\mathbb{E}[ \1_{\{\sum_{j=0}^{n-1} \RandVar_j/\sqrt{n}\ \leq
t\}}]\stackrel{\scriptscriptstyle{n\to\infty}}{\longrightarrow}
\frac{1}{\sqrt{2\pi}}\int_{-\infty}^ t e^{-x^2/2}\dd x
\end{equation}
for any $t\in\R$. To better compare these theorems, it is worth
noticing that \eqref{ascltiid} implies that \emph{almost surely}
  \begin{equation}
  \frac{1}{\log n} \sum_{k=1}^n \frac{1}{k}
  \1_{\{\sum_{j=0}^{k-1} \RandVar_j/\sqrt{k}\ \leq t\}}
  \stackrel{\scriptscriptstyle{n\to\infty}}{\longrightarrow}
  \frac{1}{\sqrt{2\pi}}\int_{-\infty}^ t e^{-x^2/2}\dd x
  \end{equation}
for any $t\in\R$. So, instead of taking the expected value, we take
a logarithmic average and obtain an almost-sure convergence.

In fact, whenever there is independence and a classical limit
theorem, the corresponding almost-sure limit theorem also holds
(under minor technical conditions), see \cite{berkes_csaki} and
references therein. The situation is more complicated for
weakly-dependent sequences, see \cite{yoshihara} and references
therein.

For dynamical systems $(X,T,m)$ given by the iteration of a map
$T:X\circlearrowleft$ which preserves the probability measure $m$,
we take $\RandVar_j=f\circ T^j$, where $f:X\to\R$ is an observable. Here,
the randomness only comes from the choice of the initial condition
$x$ according to the invariant measure of the system. The sequence
$\RandVar_j$ is identically distributed (in fact stationary), but there is
no independence in general. Nevertheless, it is well-known that many
dynamical systems display a complicated behavior which can be
adequately analyzed by probabilistic methods.

For some classes of systems, it is possible to use probabilistic
techniques for weakly dependent sequences, and prove an almost sure
invariance principle. That is, there exist $\epsilon>0$
and a Brownian motion $W$ (on a possibly extended space) such that,
almost surely,
  \begin{equation}
  \sum_{j=0}^{n-1} f\circ T^j(x) = W(n) (x)+o(n^{1/2-\epsilon}) \quad
  \text{when }n\to \infty.
  \end{equation}
This directly implies that the Birkhoff sums of $f$ satisfy an almost
sure central limit theorem, by \cite{lacey_philipp}. See
e.g.~\cite{denker:ASIP,dolgopyat:limit, melbourne:ASIP} for examples
of dynamical systems satisfying the almost sure invariance principle
-- these include Anosov maps as well as partially or non-uniformly
hyperbolic transformations.

The goal of this article is to support the idea that ``whenever we can prove a limit
theorem in the classical sense for a dynamical system, we can prove
a suitable almost-sure version''. More precisely, we will
investigate three methods that are used to prove limit theorems in
dynamical systems: spectral methods, martingale methods, and
induction arguments. We will show that whenever these methods apply,
the corresponding limit theorem admits a suitable almost-sure
version. Typically our statements will look like:
  \begin{equation}
  \frac{1}{\log n} \sum_{k=1}^n \frac{1}{k} \gdelta_{S_k f/ B_k}
  \lawto \boW \quad\textup{almost-surely}
  \end{equation}
where $f:X\to\R$ is a ``regular'' observable, $B_k$ is a suitable
normalizing sequence, $\boW$ a suitable law, and $S_k f := f+f\circ
T+\cdots +f\circ T^{k-1}$.

Let us give some more details about the methods for proving
limit theorems we mentioned above:
\begin{itemize}
\item
\textbf{Spectral methods:} If $T: X \to X$ is a probability
preserving map, the corresponding \emph{transfer operator}
is defined on $L^2$ as the adjoint of the composition by $T$. Under
suitable assumptions on the map $T$, it acts on spaces of regular
functions, and has a spectral gap. This result, which implies in
particular exponential decay of correlations, is a very useful tool
to study limit theorems for $T$. In a specific setting, the
so-called \emph{Gibbs-Markov maps}, this good spectral behavior was
used by Aaronson and Denker in \cite{aaronson_denker},
\cite{aaronson_denker:central} to prove that the suitably
renormalized Birkhoff sums of a good observable converge to a
Gaussian or stable law (see Theorem \ref{TCL_GM} for a precise
statement of their results). Under the same assumptions, we will
prove an almost sure limit theorem (Theorem \ref{thm:GM}).
It will be derived from a more general
theorem stated in terms of continuous perturbations of transfer
operators, which applies in a large variety of settings (Theorem
\ref{thm:spectral}).
\item
\textbf{Martingale methods:} Let again $T$ be a probability
preserving map on a space $X$. If $f:X\to \R$ is a function, it is
sometimes possible to write it as $f=g-g\circ T + h$, where $g$ is a
measurable function and the sequence $h\circ T^n$ is a reverse
martingale difference for some filtration. The central limit theorem
for reverse martingale differences then implies that the Birkhoff
sums of $h$ satisfy a central limit theorem. This in turn yields the
same conclusion for $f$. We will prove that a sequence of reverse
martingale differences also satisfies an almost sure central limit
theorem, by mimicking the proof in \cite{lifshits:TCLps} for the
direct martingale differences. As above, this gives an almost sure
limit theorem in the dynamical systems setting, given in Theorem
\ref{TCLps_martingales}.
\item
\textbf{Induction methods:} Let $T$ be a probability preserving map
on a space $X$, and let $Y$ be a positive measure subset of $X$. Let
$T_Y$ be the induced map on $T$, and $\phi$ the first return time.
If the Birkhoff sums of a function $f$, for the transformation
$T_Y$, satisfy a limit theorem, then it is well known (see 
e.g.~\cite{aaronson_denker_urbanski,zweimuller,melbourne_torok})
that, \emph{under suitable additional assumptions},
the function $f$
also satisfies a limit theorem for the initial map $T$. In
\cite{melbourne_torok}, the additional assumptions are formulated in
terms of the return time function $\phi$, which should essentially
satisfy a central limit theorem. Our first goal when we started to
write this paper was to extend this kind of result
to almost sure limit theorems. We were surprised to realize that
this extension was indeed possible, under \emph{weaker} assumptions.
Indeed, there is no need to assume anything on the return time
function $\phi$ (see Theorem~\ref{thm_local_global_0}).

We also tried to eliminate the conditions on $\phi$ in Melbourne and
T\"{o}r\"{o}k's classical limit theorem, and were only partially successful:
this is possible under additional assumptions on the function $f$,
which amount to a tightness condition for the maxima of the Birkhoff
sums for $T_Y$ (see Definition~\ref{def:TightMaxima}). This
condition can be checked in several practical cases, by a martingale
argument. This yields new limit theorems which could not be proved
by the previous variations around \cite{melbourne_torok}, see 
e.g.~Theorem~\ref{thm_limite_normal} in which there is no assumption on
the return time $\phi$.
\end{itemize}

We notice that almost-sure central limit theorems have been
established for certain dynamical systems in
\cite{chazottes_collet,chazottes_collet_schmitt} as a consequence of
some concentration inequalities. In fact, a strengthening of the
almost-sure central limit theorem is obtained in these papers.
Moreover, \cite{lesigne:TCLps} proves that, for \emph{any}
probability preserving dynamical system, there exists a function $f$
which satisfies the almost sure central limit theorem (but usually,
this function is quite wild).

In Section \ref{results}, we give all the precise statements of the
theorems, and the remaining sections are devoted to their proofs.

\section{Statement of main results}\label{results}

\subsection{Different notions of convergence}

In this paragraph, we modify the classical notion of convergence of
random variables in two different ways, by putting an additional
condition which will prove very useful to induce limit theorems from
a subset of the space to the whole space (see Theorem
\ref{thm_probabiliste_general}), or by studying almost sure
convergence. We will see later that these new notions of convergence
are satisfied in several cases by dynamical systems.

A continuous function $L:\R_+^*\to \R_+^*$ is \emph{slowly varying}
if, for all $\lambda>0$, $L(\lambda x)/L(x) \to 1$ when
$x\to\infty$. This implies that $L(x)=o(x^{\epsilon})$ for all
$\epsilon>0$, as well as $1/L(x)=o(x^\epsilon)$. Basic examples of
functions with slow variation are constant functions and powers of
the logarithm function.

A slowly varying function $L$ is said to be \emph{normalized} if $L$
is $C^1$ and $L'(x)=o(L(x)/x)$. Every slowly varying function is
asymptotically equivalent to a normalized slowly varying function,
see \cite[Theorem 1.3.3]{regular_variation}. In particular, if one
is only interested in the asymptotic behavior of slowly varying
functions, one can without loss of generality restrict oneself to
normalized slowly varying functions.

\begin{defn}
A \emph{renormalization function} is a function $B: \R_+^*\to
\R_+^*$ of the form $B(x)=x^d L(x)$ where $d>0$ and $L$ is a
normalized slowly varying function. The corresponding
\emph{renormalizing sequence} is $B_n:=B(n)$.
\end{defn}

\begin{defn}
\label{def:TightMaxima} Let $S_n$ be a sequence of random variables
on a probability space, and let $B_n$ be a renormalizing sequence.
We say that $(S_n/B_n, B_n)$ converges with tight maxima to a random
variable $\boW$, if $S_n/B_n$ converges in law to $\boW$, and the
sequence $M_n= (\max_{1\leq k \leq n} |S_k|)/B_n$ is tight, i.e.,
  \begin{equation}
  \forall\epsilon>0,\, \exists\ c>0 \, \textup{ s.t.}\;\forall n\geq 1,\
  \P\left\{ \max_{1\leq k \leq n} |S_k|/B_n > c\right\} \leq \epsilon\,.
  \end{equation}
\end{defn}
Notice that this property is not a property of the sequence
$S_n/B_n$ only, the renormalizing sequence $B_n$ plays a role in the
definition of $M_n$. However, abusing notations, we will usually
simply say that $S_n/B_n$ converges with tight maxima to $\boW$.

\begin{example}
\label{ex:TightMaxima} Let $\RandVar_0,\RandVar_1\dots$ be a sequence
of reverse
martingale differences. Let $S_n=\sum_{k=0}^{n-1} \RandVar_k$. Assume that,
for some renormalizing sequence $B_n$, $S_n/B_n$ converges in law to
a random variable $\boW$, and that $S_n/B_n$ is bounded in $L^1$.
Then $S_n/B_n$ also converges with tight maxima to $\boW$.
\end{example}
\begin{proof}
The maximal inequality for reverse martingales shows that, for all
$\alpha>0$ and all $n\in \N$,
  \begin{equation}
  \P\left\{ \max_{1\leq k \leq n}|S_k| \geq \alpha \right\} \leq
  \frac{C}{\alpha}\ \E(|S_n|)
  \end{equation}
where $C$ is a universal constant. In particular, for all $c>0$,
  \begin{equation}
  \P\left\{ \max_{1\leq k \leq n}|S_k|/B_n \geq c\right\} \leq
  \frac{C}{c B_n}\ \E(|S_n|)
  \end{equation}
which is bounded by $C'/c$ since $S_n/B_n$ is bounded in $L^1$.
\end{proof}

\begin{defn}
\label{def:AlmSure}
Let $S_n$ be a sequence of random variables on a
probability space, and let $B_n$ be a renormalizing sequence. We say
that $S_n/B_n$ satisfies an almost sure central limit theorem
towards a random variable $\boW$ if, for almost all $\omega$,
  \begin{equation}
  \frac{1}{\log N} \sum_{k=1}^N \frac{1}{k} \gdelta_{ S_k(\omega)/B_k}
  \lawto \boW
  \end{equation}
where $\delta_x$ is the Dirac mass at $x$, and the convergence is
the weak convergence for probability measures on $\R$.
\end{defn}
Contrary to Definition~\ref{def:TightMaxima}, this is a property of
the sequence $S_n/B_n$ only.

\subsection{Spectral arguments: Gibbs-Markov maps}\label{spectral}

\subsubsection{Almost-sure limit theorems in the i.i.d.\ case}
\label{subsubdomattr}

Let us first recall the precise statements of almost-sure limit
theorems for i.i.d.\ sequences of random variables in the domain of
attraction of a Gaussian or stable law.

A function $f$, defined on a probability space $(\Omega,\boB,m)$,
is said to \emph{belong to a domain of attraction} if 
it satisfies one the following three conditions:

\begin{enumerate}
\item[\textup{I}.] It belongs to $L^2(\Omega)$.
\item[\textup{II}.] One has $\int \1_{\{|f|>x\}} \dd m \sim
x^{-2}\ell(x)$, for some function $\ell$ such that $L(x):=2\int_1^x
\frac{\ell(u)}{u}\dd u$ is of slow variation and unbounded.
\item[\textup{III}.] There exists $p\in (1,2)$
such that $\int \1_{\{f>x\}} \dd m=(c_1+o(1))x^{-p}L(x)$
and $\int \1_{\{f<-x\}} \dd m=(c_2+o(1))x^{-p}L(x)$, where $c_1,c_2$
are nonnegative real numbers such that $c_1+c_2>0$, and $L$ is of
slow variation.
\end{enumerate}

It is convenient to say that in conditions I and II we have $p=2$,
and that $L(x)=1$ in condition I.

Let us briefly comment on these conditions. The second one is
equivalent to the fact that $\tilde{L}(x)=\int f^2 \1_{\{|f|\leq x\}}
\dd m$ is of slow variation and unbounded. Moreover, the functions
$L$ and $\tilde{L}$ then are equivalent at $+\infty$. In that case,
the function  $f$ belongs to $L^q$ for all $q<2$, but not to $L^2$.
Also $\ell(x)=o(L(x))$.

In condition III, the function $f$ belongs to $L^q$ for all $q<p$.
It may or may not belong to $L^p$, according to the behavior of the
function $L$. It never belongs to $L^q$ for $q>p$.

Note in particular that the three conditions are mutually exclusive.

The above definition of domain of attraction is motivated by the
following well-known, classical result in Probability (see 
e.g.~\cite{gnedenko_kolmogorov}):
\begin{thm}
\label{thm_limite_iid} Let $Z$ be a random variable belonging to a
domain of attraction. Let $Z_0,Z_1,\dots$ be a sequence of
independent, identically distributed, random variables with the same
law as $Z$. In all cases, we set $A_n=n\E(Z)$ and
\begin{enumerate}
\item If condition I holds, we set $B_n=\sqrt{n}$ and $\boW=
\boN(0, \E(Z^2)-\E(Z)^2)$.
\item If condition II holds, we let $B_n$ be a renormalizing sequence
with $nL(B_n)\sim B_n^2$, and $\boW=\boN(0,1)$.
\item If condition III holds, we let $B_n$ be a renormalizing sequence
such that $nL(B_n)\sim B_n^p$. Define
$c=(c_1+c_2)\Gamma(1-p)\cos\left(\frac{p\pi}{2}\right)$ and
$\beta=\frac{c_1-c_2}{c_1+c_2}$. Let $\boW$ be the law with
characteristic function
  \begin{equation}
  \E(e^{it\boW})=e^{-c|t|^p\left(1-i\beta \sgn(t)
  \tan\left(\frac{p\pi}{2}\right) \right)}.
  \end{equation}
\end{enumerate}
Then
  \begin{equation}
  \frac{\sum_{i=0}^{n-1} Z_i -A_n}{B_n} \lawto \boW.
  \end{equation}
\end{thm}

\bigskip

The conditions put on the distribution of $Z$ are almost necessary
and sufficient to get a convergence in law of that type, we only
restricted the range of $p$'s, which could also be taken in the
interval $(0,1]$.

Notice that it is possible to construct the renormalizing sequence
$B_n$ by taking $B_n=n^{1/p} \overline{L}(n)$, where
$\overline{L}$ is a normalized slowly varying
function built up from $L$, or more precisely from its
\emph{de Bruijn conjugate} (see
\cite{regular_variation}).

As a matter of fact, random variables $S_n= \sum_{i=0}^{n-1} Z_i
-A_n$ as in the statement of the preceding theorem not only converge
in law when properly rescaled, but they also converge with tight
maxima (this is a consequence of Example~\ref{ex:TightMaxima}).
Moreover,  the following theorem holds (see 
e.g.~\cite{berkes_csaki} for a proof).

\begin{thm}
\label{thm_ps_iid} Under the same hypotheses and with the same
notations as in Theorem~\ref{thm_limite_iid}, $(\sum_{i=0}^{n-1} Z_i
-A_n)/B_n$ satisfies an almost sure limit theorem towards $\boW$.
\end{thm}

\subsubsection{Almost-sure limit theorems for Gibbs-Markov maps}

In this paragraph we give the analog of Theorem~\ref{thm_ps_iid} for
Gibbs-Markov maps, which are defined as follows.

\begin{defn}
Let $T:X\circlearrowleft$ be a non-singular map on a probability,
metric space $(X,\boB,m,d)$ with bounded diameter, preserving the
probability measure $m$. This map is said to be
``\emph{Gibbs-Markov}'' if there exists a countable (measurable)
partition $\alpha$ of $X$ such that:

\begin{enumerate}
\item For all $a\in \alpha$, $T$ is injective on $a$ and $T(a)$
is a union of elements of $\alpha$.
\item There exists $\lambda>1$ such that, for all $a\in \alpha$, for
all points $x,y\in a$, $d(Tx, T y)\geq \lambda d(x,y)$.
\item Let $\Jac$ the inverse of the Jacobian of $T$. There exists
$C>0$ such that, for all $a\in \alpha$, for all points
$x,y\in a$, $\left|1-\frac{\Jac(x)}{\Jac(y)} \right| \leq C
d(Tx,Ty)$.
\item The map $T$ has the ``big image property'':
$\inf_{a\in \alpha} m(T a)>0$.
\end{enumerate}
\end{defn}

These properties say that $T$ is Markovian, uniformly expanding and
with bounded distortion. In some sense, such maps have the strongest
possible chaotic behavior, and are the first candidates when one
wants to extend a probabilistic limit theorem to dynamical systems.

Let us define the separation time, $s(x,y)$, of two points $x,y\in
X$ as the number of iterations of $T$ necessary for the orbit of $x$
and $y$ to fall into distinct atoms of the partition $\alpha$. For
$\tau<1$, define a new distance $d_\tau$ on $X$ by setting
$d_\tau(x,y)=\tau^{s(x,y)}$. If $\tau$ is sufficiently close to $1$,
the map $T$ is still Gibbs-Markov for the distance $d_\tau$.

Let $f:X \to \R$ a function. For $X' \subset X$, we let
  \begin{equation}
  Df(X')=\sup\left\{
  \frac{|f(x)-f(y)|}{d(x,y)} \ : \ x,y \in X', x\not=y\right\}\cdot
  \end{equation}
This is the best Lipschitz constant of $f$ on $X'$.

We now state a theorem asserting that the convergence results of
Theorem~\ref{thm_limite_iid} extend from the i.i.d.\ case to the
case of Gibbs-Markov maps. This result is proved in
\cite{aaronson_denker}, \cite{aaronson_denker:central} and
\cite{gouezel:stable}.

\begin{thm}
\label{TCL_GM} Let $T:X\circlearrowleft$ be a Gibbs-Markov map for a
partition $\alpha$, preserving the ergodic probability measure $m$.
Consider $f: X\to \R$ such that $\sum_{a\in \alpha} m(a)
Df(a)<\infty$ and such that the distribution of $f$ belongs to a domain of
attraction as above. Assume also $\int f \dd m=0$. Then
  \begin{equation}
  \frac{S_n f}{B_n} \lawto \boW
  \end{equation}
where $B_n=\sqrt{n}$ and $\boW=\boN(0,\sigma^2)$ for some
$\sigma^2\geq 0$ if $f\in L^2$; Otherwise $B_n$ and $\boW$ are as in
the i.i.d.\ case.
\end{thm}

We use the classical notation $S_n f:= f+f\circ T+\cdots+f\circ
T^{n-1}$.

\begin{rmq}
When $f\in L^2$, the value of $\sigma^2$ is not always $\int f^2 \dd
m$, due to the lack of independence. It is in fact equal to $\int
f^2\dd m +2 \sum_{k=1}^\infty \int f \cdot f\circ T^k \dd m$ (and
this series is converging). On the other hand, when $f\not\in L^2$,
the sequence $f\circ T^k$ behaves really as if it were independent.
\end{rmq}

In this setting, we obtain the following result concerning
almost-sure limit theorems.

\begin{thm}
\label{thm:GM} With the assumptions and notations of Theorem
\ref{TCL_GM}, $S_nf /B_n$ converges with tight maxima to $\boW$.
Moreover, it satisfies an almost sure limit theorem towards $\boW$.
\end{thm}

\subsubsection{A more general spectral result}

\label{subsub:spectral}

Theorem \ref{thm:GM} will be derived from a more general spectral
theorem which applies also to different settings. The spirit of this
paragraph is close to the ideas of \cite{hennion_herve}, with weaker continuity
assumptions.

Let $T$ be a nonsingular map on the probability space $(X,m)$ (the probability
measure $m$ is not assumed to be invariant), and let $f:X\to \R$ be
measurable.
Let $\boG$, $\boH$ be two complex Banach spaces and
let $\id: \boG \to \boH$ be a continuous linear map
such that the image of the unit ball of $\boG$ is relatively compact
in $\boH$. Assume that two elements $\alpha_0 \in \boG$ and $\ell_0
\in \boG'$ (the dual of $\boG$) are given. Finally, consider some
$0<\epsilon_0<1$, and assume that operators $\boL_t : \boG \to \boG$ are
given, for $|t| \leq \epsilon_0$.

We assume the following properties:
\begin{enumerate}
\item
For all $t,t'\in [-\epsilon_0,\epsilon_0]$ and all
$n,p\in \N$,
  \begin{equation}
  \label{assumption}
  \int e^{it S_n f \circ T^p} e^{it' S_p f} \dd m=
  \langle \ell_0, \boL_t^n \boL_{t'}^p \alpha_0 \rangle.
  \end{equation}
\item
There exist constants $C>0$, $\eta<1$ and $M\geq 1$ such that,
for all $u\in \boG$, for all $n\in \N$, for all
$t\in[-\epsilon_0,\epsilon_0]$,
  \begin{equation}
  \label{eq:LY}
  \norm{ \boL_t^n u}_{\boG} \leq C\eta^n \norm{u}_{\boG} + C M^n
  \norm{\id(u)}_{\boH}.
  \end{equation}
and
  \begin{equation}
  \norm{ \id(\boL_t^n u)}_{\boH} \leq C M^n \norm{\id(u)}_{\boH}.
  \end{equation}
\item The eigenvalues of modulus $\geq 1$ of the
operator $\boL_0$ are simple. Moreover,
$\boL_0' \ell_0=\ell_0$.
\item
  There exists $\beta_0>0$ such that, for all
$t\in[-\epsilon_0,\epsilon_0]$,
  \begin{equation}
  \norm{ \id\circ (\boL_t-\boL_0)}_{\boG \to \boH} \leq C |t|^{\beta_0}.
  \end{equation}
\end{enumerate}

It is often possible to take for $\boG$ a space of functions on
$X$. The operator $\boL_0$ is the transfer operator, $\alpha_0$ is the
function $1$ and $\ell_0$ is the integration against the measure $m$.
The perturbed operator $\boL_t$ is then
usually given
by $\boL_t (u)=\boL_0 (e^{itf}u)$
\emph{if this can be defined}. The first assumption is then a formal
consequence of the definition.
To do this, one needs to be
able to multiply an element of $\boG$ by the function $e^{itf}$, and
still get an element of $\boG$. This is not always the case. For
example, when $T$ is a Gibbs-Markov map and the function $f$ is
integrable and satisfies $\sum_{a\in\alpha} m(a) Df(a)<\infty$ (where
$\alpha$ is the Markov partition of $T$), then it is possible to
define $\boL_t$ acting on the space $\boG$
of locally H\"{o}lder functions, but not as naively as before: in
general, if $u\in \boG$, then $e^{itf} u\not\in \boG$. Nevertheless,
the operator $\boL_0$ is regularizing, and sends back $e^{itf} u$ in
$\boG$, therefore $\boL_t$ is well defined and satisfies the first assumption.

The more general setting given above is useful to treat more general
dynamical systems where the convenient spaces to act on are not spaces of
functions any more, such as in the hyperbolic setting (see
\cite{gouezel_liverani, bt:aniso}).

Notice that the second assumption is a uniform Lasota-Yorke
inequality. By Hennion's Theorem \cite{hennion_herve}, it ensures that $\boL_t$ has a
finite number of eigenvalues of modulus $\geq \rho$ for any
$\rho>\eta$, and that these eigenvalues have finite multiplicity.
The third assumption gives a more specific spectral description for
$\boL_0$.

The fourth assumption is a \emph{weak} continuity assumption. It
does not imply that $\norm{\boL_t-\boL_0}_{\boG \to \boG}\to 0$ when
$t\to 0$ (this would be a too strong assumption, which would not be
satisfied in many interesting cases, see e.g.~the case of the
stadium billiard in Paragraph \ref{subsub:stadium}). However,
together with the uniform Lasota-Yorke inequality, it is sufficient
to get continuity properties for the spectrum of $\boL_t$ by
\cite{baladi_young, keller_liverani}.

\begin{thm}
\label{thm:spectral}
Under the assumptions 1--4, let $B_n$ be a
renormalizing sequence such that $S_n f/B_n$ converges in
distribution to a
random variable $\boW$.
Then $S_n f/B_n$ satisfies an
almost sure limit theorem towards $\boW$.
\end{thm}

\subsection{Induction arguments}\label{induction}

Melbourne and T\"{o}r\"{o}k \cite{melbourne_torok} have shown that under
mild assumptions the central limit theorem for a map implies the
central limit theorem for suspension flows over that map. In fact
this holds for inducing: if an induced map satisfies a limit
theorem, so does the map on the whole space, provided the return
time is nice enough.

We will show that it is possible to replace this condition on the
return time by a condition on tight maxima. To state this result, we
need a few notations. Let $(X,\boB,m,T)$ be an ergodic dynamical
system, and let $Y\subset X$ be a subset with positive measure. For
$y\in Y$, let
  \begin{equation}
  \phi(y)=\inf\{n>0 \tq T^n y\in Y\}\ .
  \end{equation}
This is the first return time of $y$ to $Y$. For a function $f:X\to
\R$, define
  \begin{equation}
  f_Y(y)=\sum_{k=0}^{\phi(y)-1} f(T^k y)\;\;\textup{if}\;y\in Y\ ,\;
  f_Y(y)=0\;\;\textup{if}\;y\notin Y\;.
  \end{equation}
Denote by $T_Y:Y\to Y$ the induced map, that is, $T_Y y=
T^{\phi(y)}y$ for every $y\in Y$ such that $\phi(y)<\infty$. Let
$S_k^Y$ stand for the Birkhoff sums for $T_Y$. Finally set
$m_Y=m(Y)^{-1} m_{|Y}$. The map $T_Y$ is defined almost everywhere on
$Y$, and preserves the probability measure $m_Y$.

\begin{thm}\label{thm_probabiliste_general}
Let $T$ be an ergodic endomorphism of a probability space
$(X,\boB,m)$. Let $Y\subset X$ be a set with positive measure, and
let $f:X\to \R$ be an integrable function. Let
$B: \R_+^* \to \R_+^*$ be a renormalization function. Assume that
$S^Y_n f_Y/ B( n/m(Y))$ converges with tight maxima to a random
variable $\boW$, for the measure $m_Y$. 
Then $S_n f /B(n)$ converges in law to $\boW$, for the measure $m$.
\end{thm}

Under these assumptions, it is interesting to know when the
convergence of $S_n f/B_n$ to $\boW$ still has tight maxima, since
it would make it possible to induce again and again. This is the
case under a quite mild condition:
\begin{prop}
\label{Prop:InduceTM}
Under the assumptions of Theorem~\ref{thm_probabiliste_general},
assume additionally that the function $M$ defined on $Y$ by
$M(y)=\max_{1\leq k\leq \phi(y)} | S_k f(y)|$ satisfies:
  \begin{equation}
  \label{jlkshaglkja}
  \sup_{n\in \N} n\, m\{ y\in Y \tq M(y) \geq c B(n)\} \to 0 \textup{ when
  } c\to+\infty.
  \end{equation}
Then $S_n f /B(n)$ converges
with tight maxima to $\boW$.
\end{prop}
Note that the function $M$ is bounded by $|f|_Y$. For instance, if
$B(n)=\sqrt{n}$ and $|f|_Y\in L^2(Y)$, then
  \begin{equation}
  nm\{ M\geq c\sqrt{n}\}\leq
  n m\{ |f|_Y \geq c \sqrt{n}\} \leq n \E(|f|_Y^2)/(nc^2)=O(1/c^2)
  \end{equation}
which shows that the assumption \eqref{jlkshaglkja} is satisfied. More
generally, if the tails of $f_Y$ and $|f|_Y$ are comparable, then this
assumption is often satisfied.

\medskip

For the almost-sure version of those limit theorems, we will need
\emph{weaker} assumptions, since no control on the maxima will be
required:

\begin{thm}\label{thm_local_global_0}
Let $T$ be an ergodic endomorphism of a probability space
$(X,\boB,m)$. Let $Y\subset X$ be a set with positive measure, and
let $f:X\to \R$ be an integrable function. Let
$B$ be a renormalizing function. We assume that $S_n^Y f_Y/
B(n/m(Y))$ satisfies an almost sure limit theorem on $Y$, towards
$\boW$. Then $S_n f /B(n)$ also satisfies an almost sure limit
theorem towards $\boW$, on $X$.
\end{thm}

Analogues of the previous theorems hold for suspensions flows and
Poincar\'{e} sections.

\subsection{Martingale arguments}\label{martingales}

In this section we deal with the almost-sure version of the central
limit theorem due to Gordin \cite{gordin} (see also
\cite{liverani:CLT}):
\begin{thm}
\label{thm:GordinCLT}
Let $T$ be an ergodic endomorphism of a probability space
$(X,\boB,m)$. Let $\boF\subset \boB$ be a $\sigma$-algebra such that
$\boF \subset T\boF$. Consider
a square-integrable function $f: X \to \R$ such that $\int f\dd m=0$
and
  \begin{equation}
  \label{eq:HypGordin}
  \sum_{n\geq 0} \norm{ \E(f| T^n \boF) -f}_{L^2}<\infty\quad\textup{and}\quad
  \sum_{n\geq 0} \norm{ \E(f|T^{-n}\boF)}_{L^2} <\infty.
  \end{equation}
Then there exists $\sigma^2\geq 0$ such that $S_n f/\sqrt{n} \lawto
\boN(0,\sigma^2)$.
\end{thm}

We will prove in Section~\ref{section:martingales} the following
theorem.

\begin{thm}
\label{TCLps_martingales} Under the same assumptions, $S_n
f/\sqrt{n}$ converges with tight maxima to $\boN(0,\sigma^2)$.
Moreover, $S_n f/\sqrt{n}$ also satisfies an almost sure limit
theorem towards $\boN(0,\sigma^2)$.
\end{thm}

The proof of the tight maxima is essentially a rephrasing of Example
\ref{ex:TightMaxima}. On the other hand, the proof of the almost
sure limit theorem will rely on an almost-sure limit theorem for
reverse martingale differences. Since we are not aware of such a
result in the literature, we will prove it, following closely the
arguments in \cite{lifshits:TCLps} for the direct martingale
differences.

\subsection{Applications}\label{applications}

In this paragraph, we describe various dynamical systems to which
the previous results apply.

\subsubsection{Axiom A maps and flows}

Let $T:X\to X$ be the restriction of an Axiom A map to one of its
basic sets. We assume that $T$ is topologically mixing. Let $m$ be a
Gibbs measure with respect to some H\"{o}lder continuous potential.
It is well known that, if $f$ is H\"{o}lder continuous, then $S_n
f/\sqrt{n}$ converges in distribution to $\boW=\boN(0,\sigma^2)$ for
some $\sigma^2\geq 0$. Since such a transformation satisfies the ASIP,
it satisfies automatically an almost sure central limit theorem as
explained in the introduction. We nevertheless give different proofs to
show in this simple example how our theorems apply.
\begin{prop}
\label{prop:AxiomAmap}
The sequence $S_n f/\sqrt{n}$ satisfies an almost sure limit
theorem towards $\boW$.
\end{prop}
\begin{proof}
The simplest proof of the central limit theorem for $S_n f$ is
probably to show that the assumptions of Gordin's Theorem \ref{thm:GordinCLT}
are satisfied for some $\sigma$-algebra $\boF$.
This is the case if one constructs $\boF$ as follows:
fix some Markov partition of $T$, and define a set to be
$\boF$-measurable if it is
a union of local stable leaves intersected with elements of the
Markov partition.

Using this $\boF$, we can apply Theorem \ref{thm:GordinCLT} and get
the classical central limit theorem. Moreover, Theorem
\ref{TCLps_martingales} also applies, and we get the almost sure limit
theorem (as well as tight maxima).

Notice that, by using $K$-partitions as in \cite{liverani:CLT} or
\cite{dolgopyat:limit} instead of Markov partitions,
this argument extends to much more general
dynamical systems.

We could also have used Theorem \ref{TCL_GM} to prove this result,
after coding and reduction to a subshift of finite type. This argument
moreover shows that the assumption of topological mixing is not
necessary, topological transitivity would suffice.
\end{proof}

Consider now a topologically transitive Axiom A flow $T_t$ 
on a basic set $X$.
Let $m$ be a Gibbs measure with respect to
a H\"{o}lder potential. Let $f$ be a H\"{o}lder continuous function with
zero average. It is
well known that $\frac{1}{\sqrt{T}} \int_0^T f\circ T_t \dd t$
converges in distribution to a Gaussian random variable
$\boN(0,\sigma^2)$ (see e.g. \cite{melbourne_torok}).

\begin{prop}\label{hyperbolic_flows}
Under the same assumptions, for almost every $x\in X$,
  \begin{equation}
  \frac{1}{\log T} \int_1^T \dd t \ \frac{1}{t}\gdelta_{\int_0^t
  f\circ T_s(x) \dd s/\sqrt{t}}\lawto \boN(0,\sigma^2).
  \end{equation}
\end{prop}
\begin{proof}
An Axiom A flow always admits a Markov partition, and can thus be
written as a suspension over a subshift of finite type. For such a
subshift, the almost sure limit theorem is a consequence of
Proposition~\ref{prop:AxiomAmap} (or directly of
Theorem~\ref{TCLps_martingales}).
The flow version of Theorem~\ref{thm_local_global_0}
then implies the desired result for the flow.
\end{proof}

\subsubsection{Locally Gibbs-Markov maps}

Let $T:X\circlearrowleft$ be a non-singular map on a probability,
metric space $(X,\boB,m,d)$, preserving the probability measure $m$.
It is said to be \emph{locally Gibbs-Markov} if it is Markovian for
a partition $\alpha$ and if there exists $Y\subset X$ of positive
measure, which is a union of elements in $\alpha$, such that:
\begin{itemize}
\item The induced map $T_Y$ is Gibbs-Markov for the partition
$\alpha_Y= \alpha\cap Y$ and the measure $m_Y=m_{|Y}/m(Y)$.
\item For all $a\in \alpha_Y$, the return-time function $\phi$ is
constant on $a$,  equal to an integer $\phi_a \geq 1$.
\item There exists $C>0$ such that, for all $a\in \alpha_Y$, for
all $x,y\in a$, for all $0\leq k<\phi_a$, we have
$d(T^kx,T^k y)\leq C d(T^{\phi_a}x,T^{\phi_a}y)$.
\end{itemize}

In the present setting, we have the analog of Theorem
\ref{thm_ps_iid}.

\begin{thm}\label{thm_limite_normal}
Let $T:X\circlearrowleft$ be an ergodic, locally Gibbs-Markov map
for a subset $Y\subset X$. Let $f:X\to\R$ be an integrable function such that
  \begin{equation}
  \sum_{a\in \alpha} m(a)Df(a)<\infty
  \end{equation}
and
$\int f \dd m=0$. Let $f_Y:X\to \R$ be
defined for $y\in Y$ by $f_Y(y)=\sum_{k=0}^{\phi(y)-1}f(T^k y)$,
where $\phi(y)$ is the return time of $y$. If $y\not\in Y$, we set
$f_Y(y)=0$.

We assume that $f_Y$ belongs to some domain of attraction, as
defined in Paragraph~\ref{subsubdomattr}. Then $S_n f /B_n$
converges to $\boW$, and satisfies an almost sure
limit theorem towards $\boW$, where $B_n=\sqrt{n}$ and
$\boW=\boN(0,\sigma^2)$ for some $\sigma^2\geq 0$ if $f_Y\in
L^2$, and $B_n$ and $\boW$ are as in Theorem~\ref{thm_limite_iid} if
for $f_Y$ we are in the 2nd or 3rd case of that theorem.
\end{thm}
\begin{proof}
The function $f_Y$ belongs by assumption to
some domain of attraction. Moreover, it satisfies $\sum_{a\in
\alpha_Y} m_Y(a) Df_Y(a) \leq \frac{1}{m(Y)} \sum_{a\in
\alpha}m(a)Df(a)<\infty$. Hence, the assumptions of
Theorem~\ref{TCL_GM} are satisfied.
Theorem~\ref{thm:GM} then shows that $f_Y$
satisfies a limit theorem with tight maxima, and an almost sure
limit theorem. Theorems~\ref{thm_probabiliste_general}
and~\ref{thm_local_global_0} make it possible to induce these limit
theorems from $Y$ to $X$.
\end{proof}

The classical convergence result in this theorem is proved in
\cite{gouezel:skewproduct}, under suitable assumptions on $\phi$.
These assumptions can be removed here due to the notion of
convergence with tight maxima.

Young towers with summable return times, as defined in
\cite{lsyoung:annals} and \cite{lsyoung:recurrence}, are locally
Gibbs-Markov maps. More generally, several non-uniformly expanding
maps have a unique invariant absolutely continuous probability
measure and can be modelled by locally Gibbs-Markov maps. This is
for example the case for the Pomeau-Manneville maps in dimension
$1$, or the Viana maps in dimension $2$.  We refer the reader to
\cite{alves_luzzatto_pinheiro} and \cite{gouezel:viana} for more
details and general statements. Theorem~\ref{thm_limite_normal}
applies to all these examples.

\subsubsection{The stadium billiard}
\label{subsub:stadium}
The stadium billiard, or Bunimovich billiard,
has been introduced in \cite{bunimovich:stadium}. It is constituted
of two parallel segments of length $\ell$ and two semicircles of
radius $1$. The transformation is the usual billiard map in this
billiard table. It preserves the Liouville measure and is ergodic.
Let $f$ be a H\"{o}lder function with zero average. Let $I$ denote the
average of $f$ along the trajectories that bounce perpendicularly to
the segments of the billiard. It is shown in \cite{balint_gouezel}
that, if $I\not=0$, then $\frac{S_n f}{\sqrt{n\log n}}$ converges to
an explicit gaussian distribution, while if $I=0$ then $\frac{S_n
f}{\sqrt{n}}$ converges to a gaussian distribution. So, a
nonstandard normalization is needed in the first case while a
standard central limit theorem holds in the second case.

\begin{thm}
In both cases, the limit theorem admits an almost sure counterpart.
\end{thm}
\begin{proof}
In \cite{balint_gouezel}, the proof of the classical limit theorem
is given in the first case by a
spectral argument and then an induction. Using Theorems \ref{thm:spectral}
and \ref{thm_local_global_0} together with the arguments of
\cite{balint_gouezel}, we therefore obtain the desired almost sure
limit theorem.

In the second case, the proof of the classical limit theorem relies on
a martingale argument, and then on two inductions. Once again, we can
use Theorems \ref{TCLps_martingales} and \ref{thm_local_global_0} to
get the conclusion.
\end{proof}

\bigskip

The paper is organized as follows.
In Section \ref{sec:PreuveThmInducing}, we prove Theorem
\ref{thm_local_global_0}, which is the only nontrivial result of the
paper concerning convergence with tight maxima. The rest of the
paper is essentially devoted to almost sure limit theorems, with
occasional complements on convergence with tight maxima in the
different settings. More precisely, in
Section~\ref{subsec:generalites}, we establish some general results
on almost-sure limit theorems in dynamical systems that we apply
subsequently. The main result of that section, which may be of
independent interest, is an almost-sure version of a result by
Eagleson \cite{eagleson} about limit theorems to be ``mixing''. In
Section~\ref{sec:local_global}, we easily deduce
Theorem~\ref{thm_local_global_0} from the general results of
Section~\ref{subsec:generalites}. In Section~\ref{sec:spectral}, we
prove Theorem~\ref{thm:spectral}, and we show in Section~\ref{sec:GM}
how this implies the
results concerning Gibbs-Markov maps, namely
Theorem~\ref{thm:GM}. Section~\ref{section:martingales} is devoted
to the proof of Theorem~\ref{TCLps_martingales}, i.e., the almost
sure central limit theorem under Gordin's assumptions.

\subsubsection*{Acknowledgments}
We thank the anonymous referees for their very useful comments and suggestions.

\section{Inducing classical limit theorems} \label{sec:PreuveThmInducing}

In this section, we prove Theorem~\ref{thm_local_global_0},
showing that a limit theorem with tight maxima for an induced map
implies a classical limit theorem for the original map.

\begin{thm}
\label{thm:MelTor} Let $(X,\boB,m,T)$ be an ergodic probability
preserving dynamical system, and let $f:X\to \R$. Let $B_n$ be a
renormalizing sequence such that $S_n f/B_n$ converges with tight
maxima to a random variable $\boW$. Let $t_1,t_2,\dots$ be a
sequence of integer valued functions on $X$ such that $t_n/n$
converges to $1$ in probability. Let also $m'$ be a probability
measure on $X$ which is absolutely continuous with respect to $m$.
Then $S_{t_n}f/B_n$ converges in distribution to $\boW$, for the
probability measure $m'$.
\end{thm}
\begin{proof}
Fix $\epsilon>0, \delta>0$. We will show that, if $n$ is large enough,
  \begin{equation}
  \label{eq:conv0}
  m\left\{ x\tq \left|\frac{S_{t_n(x)}f(x) - S_n f(x)}{B_n} \right|
  \geq \epsilon \right\} \leq 2\delta.
  \end{equation}
This will imply that $(S_{t_n} f - S_n f)/B_n$ tends in
probability to $0$ with respect also to the measure $m'$. Since $S_n f/B_n$
converges in distribution to $\boW$, for the probability measure $m'$,
by Eagleson's Theorem \cite{eagleson},
this will conclude the proof.

Since $S_n f/B_n$ has tight maxima, there exists $c>0$ such that, for
all $n\in \N$,
  \begin{equation}
  m\left\{\max_{0\leq j\leq n} |S_j f| \geq c B_n\right\} \leq \delta.
  \end{equation}
For $z\in \R$, let $\lceil z\rceil$ denote the smallest integer $\geq
z$. Since $B_n$ is a renormalizing sequence, there exists
$\gamma\in(0,1)$ small enough that,
for all large enough $n$, $B_{\lceil 2 \gamma n \rceil} \leq
\epsilon B_n/(2c)$. We fix such a $\gamma$, and write $a_n=\lceil
(1-\gamma) n  \rceil$.

If $n$ is large enough, $m\{
x\tq |t_n(x) - n| >\gamma n\} \leq \delta$. Then
  \begin{equation*}
  m\left\{ \left| \frac{S_{t_n} f- S_n f}{B_n} \right| \geq
  \epsilon\right\}
  \leq \delta + m\left\{ \left| \frac{S_{t_n} f- S_n f}{B_n} \right| \geq
  \epsilon, t_n\in [(1-\gamma)n, (1+\gamma)n]\right\}.
  \end{equation*}
If $x$ belongs to this last set, there exists $j\in [(1-\gamma)n,
(1+\gamma)n]$ such that $|S_n f(x) - S_j f(x)|\geq \epsilon B_n$. In
particular, $| S_kf(x)-S_{a_n}f(x)|\geq \epsilon B_n/2$ for $k=j$ or
$n$. Hence,
  \begin{equation}
   m\left\{ \left| \frac{S_{t_n} f- S_n f}{B_n} \right| \geq
  \epsilon\right\}
  \leq \delta+ m\left\{ \max_{0\leq i\leq 2\gamma n}
  |S_{a_n+i}f-S_{a_n}f|\geq \epsilon B_n/2\right\}.
  \end{equation}
Since $m$ is invariant and $\epsilon B_n/2 \geq c B_{\lceil 2\gamma
n\rceil}$, the measure of this last set is at most
  \begin{equation}
  m\left\{ \max_{0\leq i \leq \lceil 2\gamma n \rceil} |S_i f| \geq c
  B_{\lceil 2\gamma n \rceil}\right\}.
  \end{equation}
This quantity is bounded by $\delta$ by definition of $c$. This
concludes the proof of \eqref{eq:conv0}.
\end{proof}

\begin{proof}[\MakeThmProofTitle{\ref{thm_probabiliste_general}}]
This result is an easy consequence of Theorem~\ref{thm:MelTor} and the
techniques of \cite{melbourne_torok} and \cite{gouezel:skewproduct},
as we will explain now. Without loss of generality, we can assume that
$T$ is invertible, since otherwise we can work in the natural
extension of $T$.

For $y\in Y$ and $N\in \N$, let $n(y,N)$ be the greatest integer $n$
such that $S_n^Y \phi(y)\leq N$. For $x\in X$, let $\pi x$ denote
its first preimage belonging to $Y$. The first two steps of the
proof of \cite[Theorem A.1]{gouezel:skewproduct} show that $S_N f(x)
/B(N)- S_{n(\pi x,N)}^Y f_Y(\pi x)/B(N)$ converges to $0$ in
probability. Hence, it is sufficient to prove that $S_{n(y,N)}^Y
f_Y(y)/ B(N)$ converges in distribution to $\boW$, for the measure
$m'$ on $Y$ with density $\dd m'=\1_Y \phi \dd m$.

By assumption, $S^Y_{\lfloor N m(Y) \rfloor} f_Y/ B(N)$ converges with
tight maxima to $\boW$, with respect to $m_Y$.
Moreover, $\int \phi \dd m_Y=1/m(Y)$ by Ka\v{c}' Formula. Hence, by Birkhoff's
ergodic Theorem, $n(y,N) \sim Nm(Y)$ for almost all $y\in Y$.
Theorem~\ref{thm:MelTor} applies and shows that $S^Y_{n(y,N)}f_Y(y)/B(N)$
converges in distribution to $\boW$ with respect to any probability
measure which is absolutely continuous with respect to $m_Y$, and in
particular for $m'$.
This concludes the proof.
\end{proof}

\begin{proof}[{\bf {\small{P}{\scriptsize ROOF
OF }{\small{P}}{\scriptsize ROPOSITION} \ref{Prop:InduceTM}}}]
For $x\in X$, let $E(x)\geq 0$ denote its first entrance time in
$Y$. Then
  \begin{multline*}
  \max_{0\leq k\leq n}|S_k f(x)|
  \\
  \leq \max_{0\leq k \leq E(x)} |S_k f(x)|
  + \max_{0\leq k \leq n} |S_k^Y f_Y (T^{E(x)}x)|
  + \max_{0\leq k \leq n} |M \circ T_Y^k (T^{E(x)}x)|.
  \end{multline*}

Let $\epsilon>0$. There exists $N\in \N$ such that $m( E(x)\geq N)
\leq \epsilon$. Therefore, for $c>0$,
  \begin{align*}
  m\left\{\max_{0\leq k\leq n} |S_k f(x)| \geq 3 c B(n)\right\}
  \leq & \, \epsilon + m\left\{\max_{0\leq k \leq E(x)} |S_k f(x)| \geq  cB(n)\right\}
  \\&
  +N m\left\{ y\in Y  \tq \max_{0\leq k \leq n} |S_k^Y f_Y (y)| \geq c B(n)\right\}
  \\&
  + N m\left\{y\in Y  \tq \max_{0\leq k \leq n} |M (T_Y^k y)| \geq c B(n)\right\}\cdot
  \end{align*}
In the upper bound, the second term is bounded by
$N \epsilon(c)$, where $\epsilon(c)$
tends to $0$ when $c\to \infty$, since $f_Y$ has
tight maxima. The last term is also bounded by $N \epsilon(c)$, by
\eqref{jlkshaglkja}. We fix $c$ so that the second and third term are
$\leq \epsilon$. Then the first term tends to $0$ when $n\to
\infty$. For large enough $n$, we get $m( \max_{0\leq k\leq n} |S_k
f(x)| \geq 3 c B(n)) \leq 4\epsilon$.
\end{proof}

\section{General results for almost-sure limit theorems in dynamics}
\label{subsec:generalites}

An almost-sure limit theorem in dynamics is a statement of the
following type: Let $T:X\circlearrowleft$ be an ergodic map preserving
a probability measure $\mum$. Let $f:X \to \R$. Under certain
assumptions, there exists a renormalizing sequence $B_n$ such that,
for almost every $x$,
  \begin{equation}
  \label{def_thmlim}
  \frac{1}{\log N} \sum_{k=1}^N \frac{1}{k} \gdelta_{S_k f(x)/B_k}
  \end{equation}
converges weakly to a probability measure on $\R$.

Let $g_n$ be a sequence of real, Lipschitz functions with compact
support which are dense (for the topology of uniform convergence) in
the space of continuous functions with compact support. The
convergence of \eqref{def_thmlim} is then equivalent to the
convergence, for each $n$, of the sequence
  \begin{equation}
  \label{def_thmlim2}
  \frac{1}{\log N} \sum_{k=1}^N \frac{1}{k}\ g_n\!\left( \frac{S_k
  f(x)}{B_k} \right)
  \end{equation}
as $N\to \infty$. For technical commodity, we will be mainly
interested in convergences like in \eqref{def_thmlim2}.

The first important observation is that the convergence in
\eqref{def_thmlim2} does not depend on the asymptotic class of
$B_k$:

\begin{lem}\label{normalisation_pasimp}
Let $x_k$ be a real sequence and let $g:\R \to \R$ be a Lipschitz function
with compact support. Assume that $\frac{1}{\log N} \sum_{k=1}^N
\frac{1}{k}g(x_k)$ converges to a limit $E$. Then, for any sequence
$\rho_k$ which tends to $1$ when $k\to \infty$,
  \begin{equation}
  \frac{1}{\log N}
  \sum_{k=1}^N \frac{1}{k}g(\rho_k x_k) \to E.
  \end{equation}
\end{lem}
\begin{proof}
It is sufficient to prove that $g(\rho_k x_k)-g(x_k) \to 0$ when
$k\to \infty$. Thus it suffices to prove that there exists a
constant $C$ such that
  \begin{equation}
  \label{supp_compact}
  \forall x\in \R,\, \forall \rho \in \R,\;
  |g(x)-g(\rho x)|\leq C |1- \rho|\ .
  \end{equation}
Let $K$ be such that $g$ is equal to zero off $[-K,K]$. If $|x|\leq
2K$, we have $|g(x)-g(\rho x)|\leq \norm{g} |x-\rho x| \leq
2K\norm{g} |1-\rho|$. If $|x|\geq 2K$ and $\rho \geq 1/2$, we have
$|\rho x| \geq |x|/2 \geq K$. Therefore, $|g(x)-g(\rho x)|=0\leq
|1-\rho|$. Finally, if $|x|\geq 2K$ and $\rho \leq 1/2$, we have
$|g(x)-g(\rho x)| \leq \norm{g}_{L^\infty} \leq 2\norm{g}_{L^\infty}
|1-\rho|$. This proves \eqref{supp_compact} in all cases.
\end{proof}

The next step is to prove that the convergence in
\eqref{def_thmlim2} is equivalent, under mild assumptions, to the
convergence of more general sums, where the normalization factor
$1/k$ is replaced by a factor of the form $\phi(T^k x)/k$. This is an
analog for almost-sure limit theorems of a result by Eagleson
\cite{eagleson}, which states that the convergence in law of $S_n f/
B_n$ for the invariant measure $\dd\mum$ is equivalent to the same
convergence for a measure $\phi\dd\mum$ with $\phi \geq 0$ and $\int
\phi \dd m=1$.

\begin{thm}
\label{ponderation_automatique} Let $T : X\circlearrowleft$ be an
ergodic map preserving a probability measure $\mum$ and let $f \in
L^1(X)$. Let $B_n$ be a renormalizing sequence. Let $g$ be a
bounded, Lipschitz function on $\R$. Then the following two
conditions are equivalent:
\begin{enumerate}
\item There exist a function $\phi \in L^1(X)$ with non-zero integral
and a set $A\subset X$ with positive measure such that, for all
$x\in A$, the quantity
  \begin{equation}
  \nu_{N,\phi,g}(x)=
  \frac{1}{\log N} \sum_{k=1}^N \frac{ \phi(T^k x)}{k}\ g\left(\frac{S_k
  f(x)}{B_k} \right)
  \end{equation}
converges to a limit $I(x)$, which may depend on $x$, when $N \to
\infty$.
\item
There exists $I\in \R$ such that, for any function $\phi \in
L^1(X)$, for almost every $x\in X$, $\nu_{N,\phi,g}(x)$ converges to
$I \int \phi\dd m$ when $N \to \infty$.
\end{enumerate}
\end{thm}

\bigskip

This theorem applies in particular when $f$ satisfies an almost-sure
limit theorem, since the first condition is then satisfied for
$\phi\equiv 1$.

The proof of this theorem relies on several technical lemmas. In the
remaining part of this section, $T$ will be an ergodic endomorphism
on a probability space $(X,\boB,\mum)$, $g$ will be a bounded,
Lipschitz function on $\R$, and $B_n$ will be a renormalizing
sequence.

\begin{lem}
\label{LemmaDeuxFonctions}
Let $\phi\in L^1(X)$ and $\psi \in
L^1(X)$. Then, for almost every $x\in X$,
  \begin{equation}
  \label{almost0}
  \frac{1}{\log N} \sum_{k=1}^N \frac{ \phi(T^k x)}{k}\min\left(1,
  \frac{|\psi(T^k x)|}{B_k}\right) \to 0.
  \end{equation}
\end{lem}
\begin{proof}
Let us first prove the lemma for $\phi=1$. Let
  \begin{equation}
  u_N(x)=\frac{1}{\log N} \sum_{k=1}^{N-1}\frac{1}{k} \min
  \left(1,\frac{ |\psi(T^k
  x)|}{B_k}\right)\cdot
  \end{equation}
We have $\int u_N \dd m= O(1/\log N)$ since $\sum 1/(k B_k)<
+\infty$. Letting $N_p = \lfloor \exp(p^2)\rfloor$, we get $\sum
\norm{u_{N_p} }_{L^1}<\infty$. Consequently, for almost every $x$,
$u_{N_p}(x) \to 0$ when $p\to \infty$. Moreover, if $N_p \leq N
<N_{p+1}$, the error made by replacing $u_{N_p}(x)$ by $u_N(x)$
tends uniformly to $0$. Hence, $u_N(x)$ tends almost everywhere to
$0$. This proves \eqref{almost0} for $\phi=1$, and consequently for
any bounded $\phi$.

If $\phi$ belongs only to $L^1$, notice that
$v_k(x)=\sum_{i=0}^{k-1} |\phi(T^i x)|$ satisfies $v_k(x)\sim
k\norm{\phi}_{L^1}$ for almost every $x$. Moreover,
  \begin{multline}
  \label{mlqksjfdlk}
  \frac{1}{\log N} \sum_{k=1}^N \frac{ |\phi(T^k x)|}{k}=
  \frac{1}{\log N} \sum_{k=1}^N \frac{
  v_{k+1}(x)-v_k(x)}{k}
  \\ = \frac{1}{\log N} \left( \frac{v_{N+1}(x)}{N} - v_1(x) +
  \sum_{k=2}^{N} v_k(x) \left( \frac{1}{k-1}-\frac{1}{k}\right)
  \right)\cdot
  \end{multline}
Hence, the limsup of this quantity is at most $\norm{\phi}_{L^1}$,
for almost every $x$.

Finally, decompose $\phi$ as $\phi_1+\phi_2$ where $\phi_1$ is
bounded and $\norm{\phi_2}_{L^1}\leq \epsilon$. Using the
convergence~\eqref{almost0} for $\phi_1$ and the previous argument
for $\phi_2$, we get that, for almost every $x$,
  \begin{equation}
  \limsup_{N\to\infty} 
  \frac{1}{\log N} \sum_{k=1}^N \frac{ \phi(T^k x)}{k}\min\left(1,
  \frac{|\psi(T^k x)|}{B_k} \right) \leq \epsilon.
  \end{equation}
Letting $\epsilon$ tend to $0$ concludes the proof.
\end{proof}

Let us note the following consequence of \eqref{mlqksjfdlk}, which
will be used several times in the sequel.

\begin{lem}
\label{lem:stableL1} If $\phi \in L^1(X)$ then for almost every
$x\in X$
  \begin{equation}
  \limsup_{N \to \infty} |\nu_{N,\phi,g}(x)| \leq \norm{g}_{L^\infty}
  \norm{\phi}_{L^1}.
  \end{equation}
\end{lem}

We now come to a more important invariance lemma.

\begin{lem}
\label{lem:phiphioT} If $\phi \in L^1(X)$ then for almost every
$x\in X$
  \begin{equation}
  \label{eq:invarT}
  \limsup_{N\to \infty} | \nu_{N,\phi,g}(x) - \nu_{N,\phi\circ T, g}(x)|=0.
  \end{equation}
\end{lem}
\begin{proof}
We first prove \eqref{eq:invarT} under the additional assumption
that $\phi$ is bounded. We have
  \begin{align*}
  \nu_{N,\phi,g}(x) - \nu_{N,\phi\circ T, g}(x)
  \!\!\!\!\!\!\!\!\!\!\!\!\!\!\!\! \!\!\!\!\!\!\!\!\!\!\!\!\!\!\!\!
   \!\!\!\!\!\!\!\!\!\!\!\!\!\!\!\! &
  \\&
  = \frac{1}{\log N}
  \left[ \sum_{k=1}^N \frac{\phi(T^k x)}{k}\ g\left( \frac{S_k
  f(x)}{B_k} \right) - \sum_{k=1}^N \frac{\phi(T^{k+1}x)}{k}\ g\left( \frac{S_k
  f(x)}{B_k} \right) \right]
  \\&
  =
  \frac{1}{\log N} \left[ \phi(Tx)\ g\left(\frac{f(x)}{B_1}\right)
  -\frac{\phi(T^{N+1}x)}{N}\ g\left( \frac{S_N f(x)}{B_N}\right)
  \right]
  \\&\ \
  +\  \frac{1}{\log N} \sum_{k=2}^N \frac{\phi(T^k x)}{k} \left[\ g\left(\frac{S_k
  f(x)}{B_k}\right) - \frac{k}{k-1}\ g\left(
  \frac{S_{k-1}f(x)}{B_{k-1}}\right) \right]\cdot
  \end{align*}
The first term tends to $0$ when $N\to \infty$, so it
suffices to estimate the second term. We have
  \begin{align*}
  \left|\ g\left(\frac{S_k
  f(x)}{B_k}\right) - \frac{k}{k-1}\ g\left(
  \frac{S_{k-1}f(x)}{B_{k-1}}\right) \right |
  \!\!\!\!\!\!\!\!\!\!\!\!\!\!\!\!\!\!\!\!\!\!\!\!\!\!\!\!\!\!
  \!\!\!\!\!\!\!\!\!\!\!\!\!\!\!\!\!\!\!\!\!\!\!\!\!\!\!\!\!\!
  \\&
  \leq
  \frac{\norm{g}_{L^\infty}}{k-1}\ +\
  \left|\ g\left(\frac{S_k
  f(x)}{B_k}\right) - g\left(
  \frac{S_{k-1}f(x)}{B_{k-1}}\right) \right |
  \\&
  \leq
  \frac{\norm{g}_{L^\infty} }{k-1}\  +\  \left|\ g\left(\frac{S_k
  f(x)}{B_k}\right) - g\left(
  \frac{S_{k}f(x)}{B_{k-1}}\right) \right |
  \\ &\quad\quad
  +
  \left|\ g\left(\frac{S_k
  f(x)}{B_{k-1}}\right) - g\left(
  \frac{S_{k-1}f(x)}{B_{k-1}}\right) \right|\cdot
  \end{align*}
We will separately estimate the contribution of each of these three
terms. First,
  \begin{equation}
  \frac{1}{\log N} \sum_{k=2}^N \frac{\phi(T^k x)}{k} \frac{\norm{g}_{L^\infty} }{k-1}
  = O(1/\log N)=o(1).
  \end{equation}
Then, for almost every $x$, $S_k f(x)=O(k)$. Hence,
  \begin{equation*}
  \left|\ g\left(\frac{S_k
  f(x)}{B_k}\right) - g\left(
  \frac{S_{k}f(x)}{B_{k-1}}\right) \right |
  \leq C |S_k f(x)| \left|\frac{1}{B_{k-1}}- \frac{1}{B_{k}}\right|
  \leq Ck \left|\frac{1}{B_{k-1}}- \frac{1}{B_{k}}\right|\cdot
  \end{equation*}
The sequence $B_k$ is eventually increasing, say, from the index $K$
on. Hence
  \begin{multline*}
  \frac{1}{\log N} \sum_{k=2}^N
  \frac{\phi(T^k x)}{k} \left|\ g\left(\frac{S_k f(x)}{B_k}\right) - g
  \left(
  \frac{S_{k}f(x)}{B_{k-1}}\right) \right |
  \\
  \leq
  \frac{1}{\log N} \sum_{k=2}^N C \left|\frac{1}{B_{k-1}}-
  \frac{1}{B_{k}}\right|
  \leq \frac{C}{\log N}\left( \sum_{k=2}^K \left| \frac{1}{B_{k-1}}-
  \frac{1}{B_k} \right| + \frac{1}{B_K}\right) =o(1).
  \end{multline*}
Finally
  \begin{equation}
  \left|\ g\left(\frac{S_k
  f(x)}{B_{k-1}}\right) - g\left(
  \frac{S_{k-1}f(x)}{B_{k-1}}\right) \right | \leq C \min\left(1,
  \frac{|f(T^{k-1}x)|}{B_{k-1}}\right)\cdot
  \end{equation}
Hence, the contribution of the corresponding term tends almost
everywhere to $0$, by Lemma~\ref{LemmaDeuxFonctions}. This concludes
the proof when $\phi$ is bounded.

To handle the case of a general $\phi\in L^1$, write
$\phi=\phi_1+\phi_2$ where $\phi_1$ is bounded and
$\norm{\phi_2}_{L^1}\leq \epsilon$. Applying the previous result to
$\phi_1$ and Lemma~\ref{lem:stableL1} to $\phi_2$, we get almost
everywhere
  \begin{multline*}
  \limsup | \nu_{N, \phi-\phi\circ T,g }(x)|
  \\
  \leq \limsup | \nu_{N, \phi_1-\phi_1\circ T,g }(x)|
  + \limsup | \nu_{N, \phi_2-\phi_2\circ T,g }(x)|
  \leq 0 + 2\epsilon\norm{g}_{L^\infty}.
  \end{multline*}
The conclusion of the lemma is obtained by letting $\epsilon$ tend
to $0$.
\end{proof}

\begin{lem}\label{lem:xTx}
If $\phi \in L^1(X)$ then for almost every $x\in X$
  \begin{equation}
  \limsup_{N\to \infty} | \nu_{N, \phi \circ T,g}(x)
  - \nu_{N,\phi,g}(Tx)|=0\ .
  \end{equation}
\end{lem}

\begin{proof}
We have
  \begin{align*}
  | \nu_{N, \phi\circ T ,g}(x) - &\nu_{N,\phi,g}(Tx)|
  \\
  & = \frac{1}{\log N} \left| \sum_{k=1}^N \frac{ \phi(T^{k+1} x)}{k}
  \left[\ g\left(\frac{S_k f(x)}{B_k} \right) -
  g\left( \frac{S_k f(Tx)}{B_k} \right) \right] \right|
  \\
  & \leq \frac{C}{\log N}
  \sum_{k=1}^N \frac{\phi(T^{k+1}x)}{k} \min\left (1, \frac{|f(T^k x) -f(x)|}{B_k}
  \right)\cdot
\end{align*}
By Lemma~\ref{LemmaDeuxFonctions}, this term converges to $0$ almost
everywhere. The lemma is proved.
\end{proof}

\begin{proof}[\MakeThmProofTitle{\ref{ponderation_automatique}}]

Let us suppose that there exists $\phi \in L^1$ whose integral is
non-zero, and such that $\nu_{N,\phi,g}(x)$ converges on a set of
positive measure. We can suppose that $\int \phi \dd m=1$.
Otherwise, replace $\phi$ by $\phi/\int \phi\dd m$. By Lemmas
\ref{lem:phiphioT} and~\ref{lem:xTx} we have, for almost every $x$,
  \begin{equation}
  \limsup_{N\to \infty} | \nu_{N,\phi,g}(x) - \nu_{N,\phi,g}(Tx)|=0\ .
  \end{equation}
In particular, the set of $x$'s for which $\nu_{N,\phi,g}(x)$
converges is invariant. Hence, by ergodicity, it is of measure one.
Moreover, the limit is an invariant function, hence a constant one.
Denote it by $I$.

Lemma~\ref{lem:phiphioT} also gives that, for all $k\in \N^*$, for
almost every $x\in X$,
  \begin{equation}
  \nu_{N ,S_k \phi/ k, g}(x) \to I.
  \end{equation}
Let $\epsilon>0$. Choose $k$ such that $\norm{ S_k \phi/ k -
1}_{L^1} \leq \epsilon$. Then, for almost every $x$,
  \begin{multline*}
  \limsup_{N\to \infty} | \nu_{N, 1, g}(x) - I|
  \\
  \leq \limsup_{N \to \infty} | \nu_{N, 1, g}(x) - \nu_{N, S_k \phi/k,
  g}(x)| + \limsup_{N \to \infty} |\nu_{N, S_k \phi/k, g}(x)-I|\ .
  \end{multline*}
The first term is at most $\epsilon\norm{g}_{L^\infty}$, by Lemma
\ref{lem:stableL1}. The second one goes to $0$. Finally, by letting
$\epsilon$ tend to $0$, we end up with: for almost every $x$,
  \begin{equation}
  \nu_{N, 1, g}(x) \to I.
  \end{equation}

Now let $\psi \in L^1(X)$ be an arbitrary function. Let
$\epsilon>0$, choose $k\in \N$ such that$\norm{S_k \psi/k - \int
\psi\dd m}_{L^1} \leq \epsilon$. Then, for almost every $x$,
  \begin{align*}
  \limsup_{N\to \infty} \left|\nu_{N,\psi, g}(x)-I\int \psi\dd m\right|
  & \leq
  \limsup_{N\to \infty} |\nu_{N, \psi, g}(x) - \nu_{N, S_k \psi/k,
  g}(x)|
  \\&\quad
  + \ \limsup_{N \to \infty} | \nu_{N, S_k \psi/k, g}(x) - \nu_{N,
  \int \psi\dd m, g}(x)|
  \\&\quad
  + \limsup_{N \to \infty} \left| \nu_{N, \int\psi\dd m,g}(x)-I\int \psi
  \dd m\right|.
  \end{align*}

The first term tends to $0$ by Lemma~\ref{lem:phiphioT}. We already
proved that the third term goes to $0$. Finally, the second one is
at most $\epsilon \norm{g}_{L^\infty}$, by Lemma~\ref{lem:stableL1}. We
conclude the proof by sending $\epsilon$ to $0$.
\end{proof}

\section{Inducing almost
sure limit theorems}\label{sec:local_global}

In this section we prove Theorem~\ref{thm_local_global_0}. For this
purpose, it suffices, according to the discussion at the beginning
of Section~\ref{subsec:generalites}, to establish the following
theorem:

\begin{thm}\label{thm_local_global}
Let $T:X\circlearrowleft$ be an ergodic map preserving a probability
measure $m$. Let $Y\subset X$ be a set of positive measure and
denote by $T_Y:Y\circlearrowleft$ the map induced by $T$, and by $\phi$ the
first return-time function. Let $f:X\to \R$ be integrable, and
define $f_Y: Y\to \R$ by
  \begin{equation}
  f_Y(y)=\sum_{k=0}^{\phi(y)-1}f(T^k y).
  \end{equation}
We will write $S_k^Y$ for the Birkhoff sums for the map $T_Y$.

Let $g:\R\to \R$ be a Lipschitz function with compact support. Let
$B$ be a renormalizing function. Assume that for almost every $y\in
Y$,
  \begin{equation}
  \label{lqskjfdlmkqsdf}
  \frac{1}{\log N} \sum_{k=1}^N \frac{1}{k}
  \ g\!\left( \frac{S^Y_{k}f_Y(y)}{B(k/m(Y))} \right) \to E
  \end{equation}
for some constant $E$. Then, for almost every $x\in X$,
  \begin{equation}
  \label{TCL_ps_sur_X}
  \frac{1}{\log N} \sum_{k=1}^N \frac{1}{k}\
  g\!\left( \frac{S_kf(x)}{B(k)} \right) \to E.
  \end{equation}
\end{thm}

\begin{proof}
By Theorem~\ref{ponderation_automatique}, it is sufficient to prove
that, for almost every $x\in Y$,
  \begin{equation}
  \label{tobeproved}
  \frac{1}{\log N} \sum_{q=1}^N \frac{1}{q}\
  g\!\left( \frac{S_q f(x)}{B(q)} \right) \to E.
  \end{equation}

Let $x\in Y$. Set $t_k(x)=\sum_{i=0}^{k-1} \phi(T_Y^i x)$: these are
the successive return times of $x$ to $Y$. For almost every $x$, we
have $t_k \sim k/m(Y)$ by Birkhoff's ergodic Theorem applied to
$T_Y$ and for the $T_Y$-invariant measure $m_Y=m_{|Y}/m(Y)$. Since
$\log t_k \sim \log t_{k+1}\sim \log k$, we are left to prove
\eqref{tobeproved} for times $N$ of the form $t_k$. We have
  \begin{equation}
  \frac{1}{\log t_k} \sum_{q=1}^{t_k} \frac{1}{q}\
  g\!\left( \frac{S_qf(x)}{B(q)} \right)
   = \frac{1}{\log t_k} \sum_{p=0}^{k-1} \sum_{q=t_{p}+1}^{t_{p+1}} \frac{1}{q}\
  g\!\left( \frac{S_q f(x)}{B(q)} \right)\cdot
  \end{equation}
For $q\in \N^*$, let $p=p(x,q)$ be the largest integer such that
$S^Y_p \phi(x)<q$. For almost all $x$, we have $S^Y_n \phi(x) \sim
n/m(Y)$, which yields $p(x,q) \sim q m(Y)$. In particular, $1/q\sim
m(Y)/p$, and $B(q)\sim B( p/m(Y))$. Using
Lemma~\ref{normalisation_pasimp}, we get
  \begin{equation*}
  \frac{1}{\log t_k} \sum_{p=0}^{k-1} \sum_{q=t_{p}+1}^{t_{p+1}} \frac{1}{q}\
  g\!\left( \frac{S_q f(x)}{B(q)} \right)
  = \frac{1}{\log k} \sum_{p=1}^{k-1} \frac{m(Y)}{p}
  \sum_{q=t_{p}+1}^{t_{p+1}} g\left( \frac{S_q f(x)}{B(p/m(Y))} \right)
  +o(1).
  \end{equation*}
For $y\in Y$, let $F(y)=\sum_{k=0}^{\phi(y)-1} |f(T^k y)|$. Then
$|S_q f(x) - S^Y_p f_Y(x)|\leq F(T_Y^p x)$. Therefore,
  \begin{multline*}
  \frac{1}{\log k} \left| \sum_{p=1}^{k-1} \frac{m(Y)}{p}
  \sum_{q=t_{p}+1}^{t_{p+1}} g\!\left( \frac{S_q f(x)}{B(p/m(Y))} \right)
  -  \sum_{p=1}^{k-1} \frac{m(Y)\phi(T_Y^p x)}{p}\
  g\!\left( \frac{S^Y_p f_Y(x)}{B(p/m(Y))} \right) \right|
  \\
  \leq \frac{C}{\log k} \sum_{p=1}^{k-1} \frac{\phi(T_Y^p x)}{p}
  \min\left(1, \frac{F(T_Y^p x)}{B(p/m(Y))}\right)\cdot
  \end{multline*}
By Lemma~\ref{LemmaDeuxFonctions}, this term tends almost everywhere
to $0$. We have proved that, for almost every $x\in Y$,
  \begin{equation*}
  \frac{1}{\log t_k} \sum_{q=1}^{t_k} \frac{1}{q}\
  g\!\left( \frac{S_qf(x)}{B(q)} \right)
  =
  \frac{1}{\log k}\sum_{p=1}^{k-1} \frac{m(Y)\phi(T_Y^p x)}{p}\
  g\!\left( \frac{S^Y_p f_Y(x)}{B(p/m(Y))} \right)
  +o(1).
  \end{equation*}

The assumption~\eqref{lqskjfdlmkqsdf} together with
Theorem~\ref{ponderation_automatique} show that this last term
converges almost everywhere to $E$.
\end{proof}

\begin{rmq}
It is possible to give a quicker proof of
Theorem~\ref{thm_local_global} by proving instead of
\eqref{tobeproved} that, for almost every $x\in Y$,
  \begin{equation}
  \label{tobeproved2}
  \frac{1}{\log N} \sum_{q=1}^N \frac{\1_Y(T^q x)}{q}\
  g\!\left( \frac{S_q f(x)}{B(q)} \right) \to Em(Y).
  \end{equation}
This is sufficient to conclude, by
Theorem~\ref{ponderation_automatique}. And there are less
computations to check~\eqref{tobeproved2} than~\eqref{tobeproved}.
However, the problem of this new proof is that it can not be
generalized easily to the case of flows, contrary to the proof given
above.
\end{rmq}

\section{Almost sure limit theorems by spectral methods}

\label{sec:spectral}

In this section, we prove Theorem \ref{thm:spectral}. Hence, we will
consider a dynamical system $(X,T,m)$ and a function $f:X\to \R$, and
assume that the assumptions 1--4 of Paragraph~\ref{subsub:spectral} hold.

\begin{lem}
\label{lem:vandermonde}
Let $(\nu_j)_{1\leq j\leq N}$ be two by two distinct complex numbers,
and let $(a_j)_{1\leq j\leq N}$ be complex numbers. Assume that
  \begin{equation}
  \sup_{n\in \N} \left| \sum_{j=1}^N a_j \nu_j^n\right| < \infty.
  \end{equation}
Then, for all $j$, either $a_j=0$ or $|\nu_j|\leq 1$.
\end{lem}
\begin{proof}
Define on the open unit disk in $\C$ an analytic function
  \begin{equation}
  \phi(z)=\sum_{n=1}^\infty \left(\sum_{j=1}^N
  a_j \nu_j^n\right) z^n.
  \end{equation}
It coincides with the function $\sum
\frac{a_j}{1-\nu_j z}$ on a small neighborhood of zero, hence on the
whole unit disk. In particular, it has no pole there. This implies
that $|\nu_j|\leq 1$ whenever $a_j\not=0$.
\end{proof}

\begin{lem}
There exist $\epsilon_1 \leq \epsilon_0$, $0<\beta\leq\beta_0$, $C>0$, $0<\rho<1$
and a function
$\lambda : [-\epsilon_1,\epsilon_1]\to\C$ such that $|\lambda(t)|\leq
1$ and, for all $t,t'\in
[-\epsilon_1,\epsilon_1]$, for all $n,p\in \N$,
  \begin{equation}
  \label{withtwoterms}
  \left| \E( e^{it S_n f\circ T^p} e^{it'S_p f}) - \lambda(t)^n
  \lambda(t')^p \right| \leq C |t|^\beta + C|t'|^\beta + C \rho^n +
  C\rho^p.
  \end{equation}
Moreover,
  \begin{equation}\label{lambdapresdeun}
  |\lambda(t)-1| \leq C |t|^\beta.
  \end{equation}
In particular,
  \begin{equation}
  \label{withoneterm}
  \left| \E( e^{it S_n f}) - \lambda(t)^n \right| \leq C|t|^\beta+C
  \rho^n.
  \end{equation}
\end{lem}
\begin{proof}
The inequality \eqref{withoneterm} (for $t'$ and $p$ instead of $t$
and $n$) is a consequence of
\eqref{withtwoterms} by taking $t=0$ and letting $n$ tend to
infinity. Hence, we just have to prove \eqref{withtwoterms}.
Estimate \eqref{lambdapresdeun} will be proved along the way.

The operator $\boL_0$ acting on $\boG$ has a simple eigenvalue at $1$,
and possibly other eigenvalues $\nu_1,\dots,\nu_k$ of modulus $\geq 1$. The
assumptions 2 and 4 of Paragraph \ref{subsub:spectral} yield the
following spectral description of $\boL_t$ for small enough $t$, by
Theorem 1 and Corollary 1 in
\cite{keller_liverani}.

The operator $\boL_t$ has an eigenvalue $\lambda(t)$ close to $1$, and
eigenvalues $\nu_j(t)$ close to $\nu_j$. Denoting by $P(t)$ and
$Q_j(t)$ the corresponding spectral projections, we can write
  \begin{equation}
  \boL_t = \lambda(t) P(t) + \sum \nu_j(t) Q_j(t) + N(t),
  \end{equation}
where $N(t)$ satisfies $\norm{N(t)^n} \leq C \rho^n$ uniformly in $t$,
for some $\rho<1$. Moreover, for any small enough $\beta>0$, we have
  \begin{equation}
  | \lambda(t)-1| \leq C |t|^\beta,\quad |\nu_j(t)-\nu_j| \leq C
|t|^\beta.
  \end{equation}
Moreover,
  \begin{equation}
  \label{borneQj1}
  \norm{ \id\circ(P(t)-P(0))}_{\boG\to \boH} \leq C|t|^\beta, \quad
  \norm{ \id \circ( Q_j(t)-Q_j(0))}_{\boG\to\boH} \leq C|t|^\beta.
  \end{equation}
Finally, the norms $\norm{P(t)}_{\boG\to\boG}$ and
$\norm{Q_j(t)}_{\boG\to\boG}$ are uniformly bounded, and
these operators satisfy
  \begin{equation}
  \label{borneQj2}
  \norm{ P(t) u}_\boG \leq C \norm{ \id( P(t)u)}_{\boH}, \quad
  \norm{ Q_j(t)u}_\boG \leq C \norm{ \id( Q_j(t)u)}_{\boH}.
  \end{equation}
To simplify the notations, we will write
$\nu_0(t)=\lambda(t)$ and $Q_0(t)=P(t)$. Let us check some algebraic
consequences of this spectral description. First, we have for any
$u\in \boG$
   \begin{equation}
  \label{oiqwureoi}
  \norm{ Q_0(0) (Q_j(t)- Q_j(0)) u}_\boG \leq C|t|^\beta.
  \end{equation}
Indeed, by \eqref{borneQj2},
  \begin{equation}
  \label{lkasfsfgdfg}
  \norm{ Q_0(0) (Q_j(t)- Q_j(0)) u}_\boG
  \leq \norm{ \id(Q_0(0) (Q_j(t)-Q_j(0)) u)}_\boH.
  \end{equation}
If $j=0$, this quantity is equal to
  \begin{equation*}
  \norm{ \id(Q_0(0) (Q_0(t)-\Id) u)}_\boH
  =\norm{ \id((Q_0(0)-Q_0(t))(Q_0(t)-\Id) u)}_\boH
  \leq C|t|^\beta,
  \end{equation*}
by \eqref{borneQj1} and the uniform boundedness of
$\norm{Q_0(t)}_{\boG\to\boG}$. On the other hand, if $j\not=0$, then
\eqref{lkasfsfgdfg} is equal to
  \begin{equation}
  \norm{ \id(Q_0(0)Q_j(t) u)}_\boH
  = \norm{ \id((Q_0(0) - Q_0(t)) Q_j(t) u)}_\boH,
  \end{equation}
which is again bounded by $C|t|^\beta$. This proves \eqref{oiqwureoi}

Let us now prove
  \begin{equation}
  \label{zxbvl;khj}
  \langle \ell_0, Q_j(t)Q_{j'}(t')\alpha_0 \rangle
  =\delta_{j0}\delta_{j'0} + O(|t|^\beta)+O(|t'|^\beta).
  \end{equation}
Since $\ell_0$ is the fixed point of $\boL_0'$, we have
  \begin{equation}
  \langle \ell_0, Q_j(t)Q_{j'}(t')\alpha_0 \rangle
  = \langle \ell_0, Q_0(0)Q_j(t)Q_{j'}(t')\alpha_0 \rangle.
  \end{equation}
Moreover,
  \begin{multline*}
  Q_0(0)Q_j(t) Q_{j'}(t')
  = Q_0(0) (Q_j(t)-Q_j(0)) Q_{j'}(t')
  \\
  + Q_0(0) Q_j(0) (Q_{j'}(t') - Q_{j'}(0))
  + Q_0(0)Q_j(0)Q_{j'}(0).
  \end{multline*}
The last term is equal to $\delta_{j0}\delta_{j'0} Q_0(0)$, while the
other ones are bounded by $C|t|^\beta$ and $C|t'|^\beta$ by
\eqref{oiqwureoi}. This proves \eqref{zxbvl;khj}.

We can now compute. We have
  \begin{equation}
  \label{expressE}
  \begin{split}
  \E( e^{it S_n f\circ T^p} e^{itS_p f})&
  = \langle \ell_0, \boL_t^n \boL_{t'}^p \alpha_0 \rangle
  \\&
  =  \sum_{j,j'=0}^k \nu_{j}(t)^n \nu_{j'}(t')^p \langle \ell_0, Q_j(t)
  Q_{j'}(t') \alpha_0 \rangle
  \\ &\ \
  + \sum_j \nu_j(t)^n \langle \ell_0, Q_j(t) N(t')^p \alpha_0 \rangle
  \\&\ \
  + \sum_{j'} \nu_{j'}(t')^p \langle \ell_0, N(t)^n Q_{j'}(t')
  \alpha_0 \rangle
  \\&\ \
  + \langle \ell_0, N(t)^n N(t')^p \alpha_0 \rangle.
  \end{split}
  \end{equation}
We will show that, in this formula, whenever there is a coefficient
$\nu_j(t)$ or $\nu_{j'}(t')$ of modulus $>1$, then the corresponding
factor vanishes. By symmetry, it suffices to do that for $\nu_j(t)$.

Fix $p\in \N$. The previous formula implies that
  \begin{equation}
  \sum_j \left(\langle \ell_0, Q_j(t) N(t')^p \alpha_0 \rangle
  + \sum_{j'} \nu_{j'}(t')^p \langle \ell_0, Q_j(t)
  Q_{j'}(t') \alpha_0 \rangle  \right) \nu_j(t)^n
  \end{equation}
is uniformly bounded, independently of $n$. Let $j$ be such that
$|\nu_j(t)|>1$. Lemma \ref{lem:vandermonde} then shows that
  \begin{equation}
  \langle \ell_0, Q_j(t) N(t')^p \alpha_0 \rangle
  + \sum_{j'} \nu_{j'}(t')^p \langle \ell_0, Q_j(t)
  Q_{j'}(t') \alpha_0 \rangle = 0.
  \end{equation}
Multiply this equation by $\rho^{-p/2}$. Then $\rho^{-p/2} \langle
\ell_0, Q_j(t) N(t')^p \alpha_0 \rangle$ is still tending to $0$, while
$\rho^{-1/2} \nu_{j'}(t')$ has modulus $>1$ for any $j'$, if $t'$ is small
enough. Applying once again Lemma \ref{lem:vandermonde} (but varying
$p$ this time), this shows
that, for all $j'$, $\langle \ell_0, Q_j(t)
Q_{j'}(t') \alpha_0 \rangle = 0$. In turn, we obtain $\langle \ell_0,
Q_j(t) N(t')^p \alpha_0 \rangle =0$. We have shown that, whenever
$|\nu_j(t)|>1$, all the corresponding factors vanish in
\eqref{expressE}.

By \eqref{zxbvl;khj}, the factor of $\lambda(t)^n \lambda(t')^p$ is
$1+O(|t|^\beta)+O(|t'|^\beta)$, which is nonzero if $t$ and $t'$ are
small enough. This yields $|\lambda(t)|\leq 1$, $|\lambda(t')|\leq
1$. 
The factors of the other terms in \eqref{expressE} are
$O(|t|^\beta)+O(|t'|^\beta)$ by \eqref{zxbvl;khj}. Hence, we have proved
  \begin{equation*}
  \E( e^{itS_n f\circ T^p} e^{it S_p f})
  = \lambda(t)^n \lambda(t')^p + O(\rho^n) + O( \rho^p) +O(|t|^\beta)
  + O(|t'|^\beta).
  \qedhere
  \end{equation*}
\end{proof}

\begin{cor}
Let $B_n \to \infty$.
The random variables $S_n f/B_n$ converge in distribution towards a
random variable $\boW$ if and only if, for all $t\in \R$,
  \begin{equation}
  \label{kl;ajsdf;lk}
  \lambda(t/B_n)^n \to \E(e^{it\boW}).
  \end{equation}
\end{cor}
\begin{proof}
The convergence in distribution of random variables is equivalent to
the pointwise convergence of the characteristic functions. That is,
$S_n f/B_n$ converges to $\boW$ if and only if, for all $t\in \R$,
  \begin{equation}
  \E(e^{it S_n f/B_n}) \to \E( e^{it \boW}).
  \end{equation}
By \eqref{withoneterm}, this is equivalent to
\eqref{kl;ajsdf;lk}.
\end{proof}
This corollary has been proved and used by Herv\'{e} in \cite{herve:kl}.

\begin{proof}[\MakeThmProofTitle{\ref{thm:spectral}}]
We will prove that, for all $t\geq 0$, for
almost all $x\in X$,
  \begin{equation}
  \label{asldjflkasfd}
  \frac{1}{\log n}\sum_{k=1}^n \frac{1}{k} \exp(it S_k f(x) /B_k) \to
  \E(e^{it \boW})
  \end{equation}
and
  \begin{equation}
  \int_{|s|\leq t}
  \frac{1}{\log n}\sum_{k=1}^n \frac{1}{k} \exp(is S_k f(x) /B_k) \dd s
  \to
  \int_{|s| \leq t} \E(e^{is \boW})\dd s.
  \end{equation}
By \cite[Lemma 6.7]{lifshits:TCLps}, this will imply the desired
almost sure limit theorem. We will only prove \eqref{asldjflkasfd},
since the other equation follows from the same estimates. So, let us
fix $t\in \R$ until the end of the proof.

We will need the following abstract lemma.
\begin{lem}
\label{AbstractLemma}
There exists a constant $C>0$ such that, for any $N\in \N$ and any
$z\in \C$, if $|z^j-1| \leq 1/2$ for all $j=0,\dots, N$, then $|z-1|
\leq \frac{C}{N}$.
\end{lem}
\begin{proof}
Write $z=re^{i\theta}$ with $r>0$ and $\theta\in [-\pi/2,\pi/2]$.
Since $r^N \in [1/2, 2]$, we get $|r-1|\leq
C/N$ for some constant $C$. If $\theta\not=0$, let $n\in \N$ be
minimal such that $|\theta n| > \pi/2$ (with $n\geq 2$ by
assumption). Then $|\theta n| \leq 2 |\theta(n-1)| \leq 2\cdot
\pi/2=\pi$. Hence, $\theta n \in [-\pi, -\pi/2) \cup (\pi/2, \pi]$. In
particular, $z^n\not\in B(1,1/2)$. This yields $n>N$. In particular,
$|\theta| \leq \pi/(2N)$.
\end{proof}

\begin{lem}
\label{ControleLambda}
There exists $C>0$ such that, for all integers $k,l$ with $l\geq Ck$, for
all $t'\in \R$ with $|t'| \leq |t|$,
  \begin{equation}
  | \lambda(t'/B_l)^k -1 | \leq C\frac{k}{l}\cdot
  \end{equation}
\end{lem}
\begin{proof}
The random variables $S_n f / B_n$ converge in distribution to
$\boW$. Hence, their characteristic functions converge, uniformly on
every compact subinterval of $\R$. In particular, there exist $N>0$
and $A>0$ such that, for all $u\in [-A,A]$ and all $n\geq N$,
  \begin{equation}
  \label{oiucxvmn}
  | \E( e^{iu S_n f/B_n} ) -1 | \leq 1/10.
  \end{equation}
Let $M$ be a large constant, and consider $l \geq 2MN$. Let $j\in
[N,l/M]$.
If $M$ is large enough, then $|t B_j /B_l|
\leq A$. Write $u=t'B_j/B_l$. Then, by \eqref{withoneterm},
  \begin{equation}
  \lambda(t'/B_l)^j = \lambda(u/B_j)^j = \E(e^{iu S_j f/B_j}) +
  O(|u|^\beta) + O( \rho^j).
  \end{equation}
Increasing $N$ and $M$ if necessary, we can ensure that the $O$ terms
are bounded by $1/10$. We get
  \begin{equation}
  \left|  \lambda(t'/B_l)^j - \E( e^{iu S_j f/B_j}) \right| \leq 1/10.
  \end{equation}
Given \eqref{oiucxvmn}, this yields
  \begin{equation}
  \label{zbmzxcvzxcv}
  \left|  \lambda(t'/B_l)^j - 1 \right| \leq 1/5.
  \end{equation}
Consider now $j\in [0,N)$. Since $N,j+N \in [N,l/M]$,
\eqref{zbmzxcvzxcv} applies to these two numbers. Therefore,
  \begin{equation*}
  \left| \lambda\left( \frac{t'}{B_l} \right)^j -1 \right|
  = \left| \frac{\bigl( \lambda(t'/B_l)^{j+N} -1 \bigr) - \bigl(
  \lambda(t'/B_l)^N -1 \bigr)}{ \lambda(t'/B_l)^N} \right|
  \leq \frac{ 1/5 + 1/5}{4/5}= \frac{1}{2}\cdot
  \end{equation*}
Hence, for all $j\in [0,l/M]$,
  \begin{equation}
  | \lambda(t'/B_l)^j -1| \leq 1/2.
  \end{equation}
By Lemma \ref{AbstractLemma}, $|\lambda(t'/B_l) - 1| \leq
C/l$. Finally, if $k\leq l/M$,
  \begin{equation*}
  \left| \lambda(t'/B_l)^k -1 \right|
  =\left| \lambda(t'/B_l) -1 \right| \left| \sum_{i=0}^{k-1}
  \lambda(t'/B_l)^i\right|
  \leq \frac{C}{l} \sum_{i=0}^{k-1} \frac{3}{2} \leq C \frac{k}{l}\cdot
  \qedhere
  \end{equation*}
\end{proof}

Let $\xi_k (x)= \exp(it S_k f(x)/B_k)-\E( \exp(it S_k f/B_k))$.
\begin{lem}
\label{LemControlXi}
For any $k, l$, we have $|\E(\xi_k \xi_l)|\leq C$. Besides, there
exist $0<\rho <1$, $\delta>0$ and $K\geq 1$
such that, if $l\geq K k$,
  \begin{equation}
  \label{ErrorTerms}
  |\E(\xi_k \xi_l)| \leq C \frac{B_k}{B_l} + C\rho^k+C
  \rho^{l-k}+\frac{C}{B_k^\delta}+ C \left( \frac{k}{l} \right)^{1/2}.
  \end{equation}
\end{lem}
\begin{proof}
Since the functions $\xi_k$ are all bounded by $2$, the first estimate
is trivial. If $k$ remains bounded, the result of the lemma is also
trivial. Hence, we may assume that $k$ is as large as needed in the
course of proof.
For the rest of the proof, we will denote by
$\epsilon(k,l)$ an error term which is compatible with
\eqref{ErrorTerms}, and we will say that such an error term is admissible.

We have
  \begin{multline*}
  \E\bigl( e^{it S_k f(x)/B_k} e^{it S_l f(x)/B_l}\bigr)
  \\
  = \E \left[ \exp\left(i\left(\frac{t}{B_k}+\frac{t}{B_l}\right) S_k
  f\right) \cdot \exp\left( i \frac{t}{B_l} S_{l-k}f \circ T^k\right)
  \right].
  \end{multline*}
Let $a=t/B_l$ and $b=t/B_k+t/B_l$. If $k$ is large enough, $a$
and $b$ are small enough so that \eqref{withtwoterms} applies. We get
  \begin{equation}
  \E\bigl( e^{it S_k f(x)/B_k} e^{it S_l f(x)/B_l}\bigr)
  = \lambda\left(\frac{t}{B_l}\right)^{l-k} \lambda\left(
  \frac{t}{B_k}+\frac{t}{B_l} \right)^k +
  \epsilon(k,l),
  \end{equation}
where $\epsilon(k,l)$ is an admissible error term.

On the other hand, by \eqref{withoneterm},
  \begin{align*}
  \E(e^{it S_k f/B_k}) \E(e^{it S_l f/B_l})
  \!\!\!\!\!\!\!\!\!\!\!\!\!\!\!\!\!\!\!\!\!\!\!\!\!\!\!\!
  \!\!\!\!\!\!\!\!\!\!\!\!\!\!\!\!\!\!\!\!
  &\\&
  = \left( \lambda\left(\frac{t}{B_k}\right)^k + O\left(
  \frac{|t|^\beta}{B_k^\beta} \right) + O(\rho^k) \right)
  \left( \lambda\left(\frac{t}{B_l}\right)^l + O\left(
  \frac{|t|^\beta}{B_l^\beta} \right) + O(\rho^l) \right)
  \\&
  =  \lambda\left(\frac{t}{B_k}\right)^k
  \lambda\left(\frac{t}{B_l}\right)^l + \epsilon(k,l).
  \end{align*}
Subtracting these two expressions, we have to show that
  \begin{equation}
  \lambda\left(\frac{t}{B_l}\right)^{l-k} \lambda\left(
  \frac{t}{B_k}+\frac{t}{B_l} \right)^k
  - \lambda\left(\frac{t}{B_k}\right)^k
  \lambda\left(\frac{t}{B_l}\right)^l
  \end{equation}
is an admissible error term to conclude. Since
$xx'-yy'=x(x'-y')+(x-y)y'$ and $\lambda$ is bounded by $1$, 
it is even sufficient to prove that
$\lambda(t/B_l)^{l-k} - \lambda(t/B_l)^l$ is admissible, as well as
$\lambda(t/B_k+t/B_l)^k - \lambda(t/B_k)^k$.

Let us first study $\lambda(t/B_l)^{l-k} - \lambda(t/B_l)^l$.
Since $\lambda(t/B_l)^l$ is uniformly bounded, it suffices
to prove that $\lambda(t/B_l)^{-k}-1$ is admissible. It even suffices
to prove that $\lambda(t/B_l)^{k}-1$ is admissible.
This is a consequence of Lemma \ref{ControleLambda} if $l \geq K k$
for some large enough $K$.

Let us now turn to $\lambda(t/B_k+t/B_l)^k -
\lambda(t/B_k)^k$. We have
  \begin{align*}
  \left| \lambda\left(\frac{t}{B_k}+\frac{t}{B_l}\right)^k -
  \lambda\left(\frac{t}{B_k}\right)^k\right|
  \!\!\!\!\!\!\!\!\!\!\!\!\!\!\!\!\!\!\!\!\!\!\!\!\!\!\!
  &\\&
  = \left|
  \E\left( e^{i\left(\frac{t}{B_k}+\frac{t}{B_l}
  \right) S_k f} \right)-\E\left( e^{i\frac{t}{B_k}
  S_k f} \right) \right|+ \epsilon(k,l)
  \\&
  \leq \E\left| e^{i\frac{t}{B_l} S_k f}-1\right| + \epsilon(k,l).
  \end{align*}
Moreover,
  \begin{align*}
  \left(\E\left| e^{i\frac{t}{B_l} S_k f}-1\right|\right)^2
  &
  \leq \E \left( \left| e^{i\frac{t}{B_l} S_k f}-1\right|^2 \right)
  = \E \left( 1-e^{i\frac{t}{B_l} S_k f} + 1- e^{-i\frac{t}{B_l} S_k
  f}\right)
  \\&
  = 1- \lambda(t/B_l)^k + 1-\lambda(-t/B_l)^k + \epsilon(k,l).
  \end{align*}
By Lemma \ref{ControleLambda}, $\left|1- \lambda(t/B_l)^k\right| \leq C
k/l$ and $ \left|1- \lambda(-t/B_l)^k\right| \leq C
k/l$.
 This concludes the proof of
the lemma.
\end{proof}

Let
  \begin{equation}
  T_N =\frac{1}{\log N} \sum_{k=1}^N \frac{1}{k} \left(e^{it S_k
  f/B_k} - \E(e^{itS_k f/B_k}) \right).
  \end{equation}
Then $T_N =\frac{1}{\log N} \sum \frac{\xi_k}{k}$. Using Lemma
\ref{LemControlXi}, and the estimates $\sum_{l=k}^\infty
\frac{1}{lB_l} \sim \frac{C}{B_k}$ and $\sum_{l=k}^\infty
\frac{1}{l\sqrt{l}} \sim \frac{C}{\sqrt{k}}$ 
(see \cite[Proposition 1.5.10]{regular_variation}), we get
  \begin{align*}
  \E(|&T_N|^2)
  \\&
  \leq \frac{C}{(\log N)^2} \sum_{k=1}^N \frac{1}{k} \left(
  \sum_{k\leq l \leq Kk} \frac{1}{l} +
  \sum_{l=Kk}^N \frac{1}{l} \left( \frac{B_k}{B_l} + \rho^k + \rho^{l-k}
  + \frac{1}{B_k^\delta}+ \frac{\sqrt{k}}{\sqrt{l}}\right) \right)
  \\&
  \leq \frac{C}{(\log N)^2} \sum_{k=1}^N \frac{1}{k} \left(\log K+
  \frac{B_k}{B_k} + \rho^k \log N + C + \frac{1}{B_k^\delta} \log N+
  \frac{\sqrt{k}}{\sqrt{k}}\right).
  \end{align*}
Since $\frac{1}{kB_k^\delta}$ is summable, as well as
$\frac{\rho^k}{k}$, we obtain
  \begin{equation}
  \E( |T_N|^2) \leq \frac{C}{\log N}\cdot
  \end{equation}
Let $N_p = \lfloor \exp(p^2) \rfloor$. Since $\E( |T_{N_p}|^2)$ is
summable, $T_{N_p}(x)$ converges to $0$ almost everywhere. That is,
for almost all $x$,
  \begin{equation}
  \frac{1}{\log N_p} \sum_{k=1}^{N_p} \frac{1}{k} \left(e^{it S_k
  f(x)/B_k}\right) \to \E(e^{it\boW}).
  \end{equation}
For a general $N$, we choose $p$ such that $N_p \leq N<N_{p+1}$ and
check that the difference between the previous sums for $N$ and for
$N_p$ converges to $0$. This concludes the proof of Theorem
\ref{thm:spectral}.
\end{proof}

\section{Limit theorems for Gibbs-Markov maps}\label{sec:GM}

In this section, we prove Theorem~\ref{thm:GM}, by applying
Theorem~\ref{thm:spectral}. Let us first recall some useful facts on
Gibbs-Markov maps.

Denote by $\boG$ the set of bounded, locally Lipschitz functions $u$
(i.e., functions $u$ satisfying $\sup_{a\in\alpha}Du(a)<\infty$ and
$\norm{u}_{L^\infty}<\infty$), endowed with its canonical norm
  \begin{equation}
  \norm{u}_{\boG}=\sup_{a\in\alpha}Du(a) + \norm{u}_{L^\infty}.
  \end{equation}
The transfer operator $\boL$ associated to $T$ acts on $\boG$ and
satisfies a Lasota-Yorke inequality
  \begin{equation}
  \norm{ \boL^n u}_\boG \leq C \eta^n \norm{u}_\boG+ C
  \norm{u}_{L^1}.
  \end{equation}
Let $f$ be a function satisfying the assumptions of
Theorem~\ref{thm:GM}. By \cite{aaronson_denker,
aaronson_denker:central, gouezel:stable}, it is possible to define an
operator $\boL_t$ acting on $\boG$ by $\boL_t(u)=\boL(
e^{itf}u)$. Moreover, it satisfies the assumptions 1--4 of
Paragraph~\ref{subsub:spectral} for $\boH=L^1$ and $\id$ the canonical 
inclusion. Hence, Theorems~\ref{TCL_GM} and
\ref{thm:spectral} together show that $S_n f/B_n$ satisfies an almost
sure central limit theorem.

To prove the tight maxima statement of Theorem~\ref{thm:spectral}, we
will need a more precise description of the mixing properties of $\boL$.
Let $r$ be the $\gcd$ of the return times of an atom of the partition
into itself. If $r=1$, the map $T$ is mixing and its correlations
decrease exponentially fast: for every function $u\in \boG$,
$\norm{\boL^n u -\int u\dd m}_\boG \leq C\eta^n \norm{u}_\boG$, for some
$\eta<1$. When $r>1$, there exists a partition of $X$, say
$X_0,\dots, X_{r-1}$, each $X_i$ being a union of elements of
$\alpha$, such that $T$ maps $X_i$ to $X_{i+1}$, for every $i\in
\Z/r\Z$. Let $\boB_0$ be the (finite) $\sigma$-algebra generated by
$\{X_0,\dots,X_{r-1}\}$, and $\Pi: u\mapsto \E(u| \boB_0)$. The
operator $\Pi$ is the projector on the eigenvalues with modulus $1$
of $\boL$. In particular, for every function $u\in \boG$,
  \begin{equation}
  \label{eq:SpectralGap}
  \bigl\|\boL^n (u -\Pi u)\bigr\|_{\boG} \leq C \eta^n
  \norm{u}_{\boG}
  \end{equation}
for some $\eta<1$. Since $\Pi$ is a conditional-expectation
operator, it satisfies $\int \Pi(u) v \dd m=\int \Pi(u)\Pi(v) \dd m$. Lastly,
$\Pi$ and $\boL$ commute.

If $u$ is integrable and $\sum m(a)Du(a)<\infty$, then
  \begin{equation}
  \label{envoie_dans_boL}
  \boL u \in \boG, \text{ and }\bigl\|\boL u\bigr\|_\boG\leq C\left(
  \sum_{a\in\alpha} m(a)Du(a)+ \int |u|\dd m\right)\, .
  \end{equation}

\begin{lem}\label{lem:borne_vite}
Assume that $f$ belongs to some domain of attraction with $\int f\dd m=0$.
Let $B_n$ be
the renormalizing sequence given by Theorem~\ref{thm_limite_iid} for
this domain of attraction. Then $S_n f/B_n$ is bounded in $L^1$.
\end{lem}
\begin{proof}
The result is clear for $\Pi f$, since its Birkhoff sums are
bounded. So, without loss of generality, we can replace $f$ by
$f-\Pi f$ and suppose that $\Pi f=0$.

We will use the following estimates on the probabilities. They are
the consequence of the slow variation of $L$ and of the choice of
$B_n$. They are easy to verify for the three types of domain of
attraction defined in Paragraph~\ref{subsubdomattr}.
  \begin{equation}
  \label{recapitule_estimees}
  \int \1_{\{|f|>B_n\}}\dd m \leq \frac{C}{n}\;;
  \int |f| \1_{\{|f|>B_n\}}\dd m \leq \frac{C B_n}{n}\;;
  \int f^2 \1_{\{|f|\leq B_n\}}\dd m  \leq \frac{C B_n^2}{n}\cdot
  \end{equation}

We define a function $\phi_n$ by $\phi_n(x)=f(x)$ if $|f(x)|\leq
B_n$, and $\sgn(f(x))B_n$ otherwise. It satisfies $D \phi_n(a) \leq
Df(a)$. Besides, since $\int f\dd m=0$, \eqref{recapitule_estimees}
implies that $\left| \int \phi_n\dd m\right| \leq \frac{C B_n}{n}$. Then
let $\psi_n=\phi_n-\Pi \phi_n$ and $\chi_n=f-\psi_n$. Since $\Pi
f=0$, they satisfy $\norm{\chi_n}_{L^1}=O(B_n/n)$ and
  \begin{equation}
 \label{borne_L2}
  \norm{\psi_n}_{L^2}^2 \leq \int \phi_n^2 \dd m
  \leq \int f^2 \1_{\{|f|\leq B_n\}}\dd m + B_n^2 \int
  \1_{\{|f|>B_n\}} \dd m\leq C \frac{B_n^2}{n}\cdot
  \end{equation}

We have $\int |S_n f|\dd m \leq \int |S_n \psi_n|\dd m + \int |S_n \chi_n|\dd m$.
Since $\norm{\chi_n}_{L^1}=O(B_n/n)$, we have
  \begin{equation}
  \label{eq:borne_psi2}
  \int |S_n \chi_n|\dd m \leq C B_n.
  \end{equation}
Moreover,
  \begin{equation}
  \int |S_n \psi_n|^2 \dd m=  n \int (\psi_n)^2 \dd m+ 2 \sum_{k=1}^n
  (n-k) \int \psi_n \cdot \psi_n \circ T^k \dd m.
  \end{equation}
The first term is bounded by \eqref{borne_L2}. For the second one,
  \begin{equation}
  \left|\int \psi_n\cdot \psi_n \circ T^k \dd m\right|
  =\left|\int \boL^k \psi_n \cdot \psi_n \dd m\right|
  \leq C \eta^{k-1} \norm{ \boL \psi_n}_{\boG}
  \end{equation}
by \eqref{eq:SpectralGap} and because $\Pi\psi_n=0$.
The $L^1$ norm of $\psi_n$ is bounded, as well as $\sum_{a\in\alpha}
m(a) D\psi_n(a)$. Hence \eqref{envoie_dans_boL} shows that $\boL
\psi_n$ is uniformly bounded in $\boG$. By gathering the preceding
equations we end up with
  \begin{equation}
  \label{eq:borne_psi1_interm}
  \int |S_n \psi_n|^2 \dd m\leq C n \frac{B_n^2}{n}+C\sum_{k=1}^n
  (n-k)\eta^k \leq C B_n^2.
  \end{equation}
Equations \eqref{eq:borne_psi2} and \eqref{eq:borne_psi1_interm} show that
$S_n f / B_n$ is bounded in $L^1$. This proves the lemma.
\end{proof}

\begin{proof}[\MakeThmProofTitle{\ref{thm:GM}}]
We have already proved the almost sure central limit theorem, it only
remains to check the tight maxima statement.
This statement is trivial for $\Pi f$, since its
Birkhoff sums are uniformly bounded. Hence, we can without loss of
generality replace $f$ with $f-\Pi f$, and assume that $\Pi f=0$.

Let $g=\sum_{n=1}^\infty \boL^n f$. This series is convergent in
$\boG$ by the spectral gap property \eqref{eq:SpectralGap}, and
\eqref{envoie_dans_boL}. The
function $h=f+g-g\circ T$ then satisfies $\boL h=0$, i.e., $\E( h|
T^{-1}\boB)=0$. The sequence $h\circ T^n$ is therefore a reverse martingale
difference for the filtration $\boF_n=T^{-n}\boB$. Moreover, $S_n
h/B_n=S_n f/B_n + (g-g\circ T^n)/B_n$. By Lemma \ref{lem:borne_vite},
$S_n f/B_n$ is bounded in $L^1$, hence $S_n h/B_n$ is also bounded in
$L^1$. Consequently,
Example \ref{ex:TightMaxima} shows that $S_n h/B_n$ has tight
maxima.
To conclude, we have to show that the sequence $(g-g\circ
T^n)/B_n$ also has tight maxima. This is a consequence of the
boundedness of $g$.
\end{proof}

\section{Almost sure limit theorems by martingale arguments}

\label{section:martingales}

\subsection{Almost-sure limit theorem for reverse martingale differences}

\begin{thm}
\label{TCLps_martingale} Let $\boF_n$ be a decreasing sequence of
$\sigma$-algebras on a probability space, and let $\RandVar_n$ be a
$\boF_n$-measurable square-integrable random variable such that
$\E(\RandVar_n | \boF_{n+1})=0$. Let $B_n\in \R_+$ increase to infinity,
let $\zeta$ be a non-negative random variable, and let $b_k\in \R_+$
be a bounded sequence with $\sum b_k=+\infty$. Assume that
\begin{enumerate}
\item Almost surely, $\RandVar_n/ B_n \to 0$.
\item Almost surely, $\frac{1}{\sum_{k=1}^n b_k} \sum_{k=1}^n b_k
\gdelta_{ \sum_{j=1}^k \RandVar_j^2(\omega) / B_k^2}$ converges weakly to
$\delta_{\zeta(\omega)}$.
\item The sequence $b_k$ satisfies
$b_k=O\left( \frac{B_k-B_{k-1}}{B_k}\right)\cdot$
\item We have
  \begin{equation}
  \sup_{k} \E\bigl[ \max_{1\leq j \leq k} |\RandVar_j|^2 / B_k^2 \bigr]
  <\infty.
  \end{equation}
\end{enumerate}
Then, for almost all $\omega$, the real measure
  \begin{equation}
  \frac{1}{\sum_{k=1}^n b_k} \sum_{k=1}^n b_k \gdelta_{ \sum_{j=1}^k
  \RandVar_j(\omega) / B_k}
  \end{equation}
converges weakly to the measure $\boN(0, \zeta(\omega))$.
\end{thm}
The proof will follow closely \cite{lifshits:TCLps}, except that we
deal with reverse martingales instead of martingales, which means we
have to reverse all the stopping time arguments.

\begin{proof}
Let $M>1$. Define a stopping time
  \begin{equation}
  \tau_k = \max\bigl[1, \sup \{ 1\leq l \leq k \tq \sum_{j=l}^{k} \RandVar_j^2 >
2M B_k^2\}\bigr].
  \end{equation}
The set $\{ \tau_k=j\}$ is $\boF_j$-measurable. Let now
  \begin{equation}
  \RandVar'_{jk}= \RandVar_j \1_{\{\tau_k \leq j\}}.
  \end{equation}
Since $\1_{\{\tau_k \leq j\}}$ is $\boF_{j+1}$-measurable, $\E(\RandVar'_{jk} |
\boF_{j+1})=0$. We will prove that, for almost all $\omega$ with
$\zeta(\omega)\leq M$,
  \begin{equation}
  \label{eq:toprove}
  \frac{1}{\sum_{k=1}^n b_k}\sum_{k=1}^n b_k \gdelta_{\sum_{j=1}^k
  \RandVar'_{jk}/B_k} \to \boN(0,\zeta(\omega)).
  \end{equation}

Let us show that this convergence implies the theorem. The
difference between $\frac{1}{\sum_{k=1}^n b_k}\sum_{k=1}^n b_k
\gdelta_{\sum_{j=1}^k  \RandVar'_{jk}/B_k}$ and $\frac{1}{\sum_{k=1}^n
b_k}\sum_{k=1}^n b_k \gdelta_{\sum_{j=1}^k
  \RandVar_j/B_k}$ has total mass at most
  \begin{equation}
  \frac{2}{\sum_{k=1}^n b_k} \sum_{k=1}^n b_k \1_{\{\tau_k>1\}}
  \leq \frac{2}{\sum_{k=1}^n b_k} \sum_{k=1}^n b_k \1_{\{\sum_{j=1}^k
  \RandVar_j^2 > 2M B_k^2\}}.
  \end{equation}
Let $f:\R\to \R$ be the piecewise affine function equal to $0$ on
$\{ x\leq M\}$ and to $1$ on $\{x \geq M+1\}$. Then this total mass
is at most
  \begin{equation}
  2 \int_\R f(x) \dd\left[ \frac{1}{\sum_{k=1}^n b_k} \sum_{k=1}^n b_k
  \gdelta_{ \sum_{j=1}^k \RandVar_j^2(\omega) / B_k^2}\right] (x).
  \end{equation}
For almost all $\omega$, this quantity converges by assumption to
$2\int_\R f(x) \dd[ \gdelta_{\zeta(\omega)} ](x)$, which is zero when
$\zeta(\omega)\leq M$. Hence, the conclusion of the theorem holds for
almost all $\omega$ with $\zeta(\omega)\leq M$. Taking a sequence
$M_n=n$ and since $\zeta$ is finite almost surely, we obtain the full
conclusion of the theorem.

So, we just have to prove \eqref{eq:toprove}. We will rather prove
that, for all $t\geq 0$, for almost all $\omega$ with
$\zeta(\omega)\leq M$,
  \begin{equation}
  \label{eq:toprove1}
  \frac{1}{\sum_{k=1}^n b_k}\sum_{k=1}^n b_k \exp\left(it \sum_{j=1}^k
  \RandVar'_{jk}(\omega)/B_k\right) \to \exp(-\zeta(\omega) t^2/2)
  \end{equation}
and
  \begin{equation}
  \label{eq:toprove2}
  \int_{|s|\leq t}
  \frac{1}{\sum_{k=1}^n b_k}\sum_{k=1}^n b_k \exp\left(is \sum_{j=1}^k
  \RandVar'_{jk}(\omega)/B_k\right) \dd s \to \int_{|s|\leq t} \exp(-\zeta(\omega)
  s^2/2)\dd s.
  \end{equation}
By \cite[Lemma 6.7]{lifshits:TCLps}, this will imply the desired
convergence \eqref{eq:toprove} almost surely. In fact, we will only
prove \eqref{eq:toprove1}, since \eqref{eq:toprove2} follows from
the same estimates.

There exists a function $r:\R \to \C$ with $|r(x)|\leq C |x|^3$ such
that
  \begin{equation}
  \exp(itx)=\exp(-t^2x^2/2) (1+itx) \exp(r(x)).
  \end{equation}
We obtain
  \begin{multline*}
  \exp\left(it \sum_{j=1}^k \RandVar'_{jk}/B_k\right)
  \\ =
  \exp\left( -t^2 \sum_{j=1}^k {\RandVar'_{jk}}^{\!\!2}/2B_k^2
  \right) \prod_{j=1}^k (1+it \RandVar'_{jk}/B_k) \exp\left( \sum_{j=1}^k
  r(t\RandVar'_{jk}/B_k) \right).
  \end{multline*}
We will denote this last product by $\boE_k(t) \Pi_k(t) R_k(t)$.
Writing
  \begin{multline*}
  \boE_k(t) \Pi_k(t) R_k(t)-\exp(-\zeta(\omega) t^2/2)=
  (\boE_k(t)-\exp(-\zeta(\omega) t^2/2))\Pi_k(t)R_k(t) \\
  + \exp(-\zeta(\omega) t^2/2) \Pi_k(t) (R_k(t)-1) + \exp(-\zeta(\omega) t^2/2) (\Pi_k(t)-1),
  \end{multline*}
we get
  \begin{equation}
  \label{decomposition_en_trois}
  \begin{split}
  \Biggl|\frac{1}{\sum_{k=1}^n b_k}\sum_{k=1}^n b_k& \exp\left(it
  \sum_{j=1}^k
  \RandVar'_{jk}(\omega)/B_k\right) - \exp(-\zeta(\omega) t^2/2)\Biggr|
  \\&
  \leq  \frac{1}{\sum_{k=1}^n b_k}\sum_{k=1}^n b_k
  |\boE_k(t)-\exp(-\zeta(\omega) t^2/2)| |\Pi_k(t)R_k(t)|
  \\&\quad
  +\frac{1}{\sum_{k=1}^n b_k}\sum_{k=1}^n b_k \exp(-\zeta(\omega) t^2/2)
  |\Pi_k(t) (R_k(t)-1)|
  \\&\quad
  + \frac{1}{\sum_{k=1}^n b_k}\left|\sum_{k=1}^n b_k \exp(-\zeta(\omega) t^2/2)
  (\Pi_k(t)-1)\right|\ .
  \end{split}
  \end{equation}
If we can prove that these three terms tend to $0$ for almost all
$\omega$ with $\zeta(\omega)\leq M$, we will have proved
\eqref{eq:toprove1} and the proof will be complete.

Write $N_k=\max_{1\leq j\leq k}|\RandVar_j|$. We have
  \begin{align*}
  \left| \sum_{j=1}^k
  r(t\RandVar'_{jk}/B_k) \right| &
  \leq \sum_{j=1}^k C t^3 |\RandVar'_{jk}|^3/B_k^3
  \leq Ct^3 \left( \RandVar_{\tau(k)}^2+ \sum_{j=\tau(k)+1}^k \RandVar_j^2 \right)
  N_k /B_k^3
  \\ &
  \leq Ct^3 \left(N_k^2 +2M B_k^2\right) N_k/B_k^3
  =C t^3 \left(N_k^2/B_k^2+ 2M\right) N_k/B_k.
  \end{align*}
For almost all $\omega$, $N_k/B_k \to 0$ by assumption. Hence,
almost surely, $R_k(t)$ tends to $1$ (and is in particular bounded).
In the same way,
\begin{align*}
|\Pi_k(t)|^2= \prod_{j=1}^k (1 +t^2 {\RandVar'_{jk}}^{\!\!2}/B_k^2) & \leq
\exp\left( t^2\left(\RandVar_{\tau(k)}^2 + \sum_{j=\tau(k)+1}^k
  \RandVar_j^2\right)/B_k^2 \right)\\
& \leq \exp( t^2(N_k^2/B_k^2 + 2M))\ .
\end{align*}
Consequently, $\Pi_k(t)$ is almost surely bounded. This proves that
the second term in \eqref{decomposition_en_trois} tends almost
surely to $0$. Moreover, almost surely, the first term in
\eqref{decomposition_en_trois} is bounded by
$\frac{C}{\sum_{k=1}^nb_k} \sum_{k=1}^n b_k |\boE_k(t)
-\exp(-\zeta(\omega) t^2/2)|$, which is at most
  \begin{multline*}
  \frac{C}{\sum_{k=1}^nb_k} \sum_{k=1}^n b_k \left[ \1_{\{\tau_k>1\}} +
  \left|\exp( -t^2\sum_{j=1}^k \RandVar_j^2/2B_k^2)- \exp(-\zeta(\omega)t^2/2)
  \right|\right]=
  \\
  \frac{C}{\sum_{k=1}^n b_k}  \sum_{k=1}^n b_k \1_{\{\tau_k>1\}} \,\, + \\
  C \int_\R |\exp(-x t^2/2) -\exp(-\zeta(\omega)t^2/2)| \dd\left[
  \frac{1}{\sum_{k=1}^n b_k} \sum_{k=1}^n b_k \gdelta_{\sum_{j=1}^k
  \RandVar_j^2/B_k^2} \right](x).
  \end{multline*}
We have seen that the first term tends to $0$ for almost all
$\omega$ such that $\zeta(\omega)\leq M$. Moreover, the second term
converges almost surely to $\int_\R |\exp(-x t^2/2)
-\exp(-\zeta(\omega)t^2/2)| \dd \gdelta_{\zeta(\omega)} (x)=0$. This
proves that the first term in \eqref{decomposition_en_trois} tends
almost surely to $0$ on $\{ \zeta \leq M\}$.

The third term in \eqref{decomposition_en_trois} is more complicated
to deal with. Notice that $\E(\Pi_k(t)| \boF_2)=\prod_{j=2}^k (1+it
\RandVar'_{jk}/B_k) \E(1+\RandVar'_{1k} | \boF_2)= \prod_{j=2}^k (1+it
\RandVar'_{jk}/B_k)$ since $\E(\RandVar'_{jk} | \boF_{j+1})=0$.
By induction, we
get $\E( \Pi_k(t) | \boF_m) = \prod_{j=m}^k (1+it \RandVar'_{jk}/B_k)$. In
particular,
  \begin{equation}
  \E(\Pi_k(t))=1.
  \end{equation}
For $l\leq k$, let us estimate $\E(\Pi_k(t) \overline{\Pi_l(t)})$.
For $p\geq 1$, write
  \begin{equation}
  A_p= \prod_{j=p}^l (1-it \RandVar'_{jl}/B_l) \prod_{j=p}^k (1+it
  \RandVar'_{jk}/B_k).
  \end{equation}
For $p>l$, there is no $\RandVar'_{jl}$ term. In particular, the same
argument as above shows that $\E(A_p)=1$. Consider now $p\leq l$.
Then
  \begin{align*}
  \E(A_p | \boF_{p+1})&= A_{p+1} \E\bigl( (1-it\RandVar'_{pl}/B_l)(1+it \RandVar'_{pk}/B_k)
  | \boF_{p+1} \bigr)
  \\&
  =A_{p+1} \Biggl[ 1 -it\E(\RandVar'_{pl}/B_l | \boF_{p+1}) +it \E(\RandVar'_{pk}/B_k
  | \boF_{p+1})
  \\ &\quad \quad \quad \quad+\frac{t^2}{B_k B_l} \E(\RandVar'_{pk} \RandVar'_{pl} |
  \boF_{p+1})\Biggr]
  \\&
  =A_{p+1} + A_{p+1} \frac{t^2}{B_k B_l} \E(\RandVar'_{pk} \RandVar'_{pl} |
  \boF_{p+1}).
  \end{align*}
Taking expectations, we get
  \begin{equation}
  | \E(A_{p+1})-\E(A_p)| = \frac{t^2}{B_kB_l}
  \Bigl| \E(A_{p+1} \1_{\{\tau_k \leq p\}} \1_{\{\tau_l
  \leq p\}} \E(\RandVar_p^2 | \boF_{p+1} )) \Bigr|.
  \end{equation}
If $\tau_k\leq p$ and $\tau_l\leq p$, we have
  \begin{align*}
  |A_{p+1}|^2 &\leq  \prod_{j=\tau_k+1}^k (1+t^2 \RandVar_j^2/ B_k^2)
  \prod_{j=\tau_l+1}^l (1+t^2 \RandVar_j^2/B_l^2) \\
  & \leq \exp\left( t^2 \sum_{j=\tau_k+1}^k  \RandVar_j^2/ B_k^2\right)
  \exp\left( t^2 \sum_{j=\tau_l+1}^l \RandVar_j^2/B_l^2 \right)
  \\
  & \leq  \exp( 2Mt^2) \exp(2Mt^2).
  \end{align*}
Hence,
  \begin{equation}
  | \E(A_{p+1})-\E(A_p)| \leq   \frac{t^2 \exp(2Mt^2)}{B_kB_l}
  \ \E(\RandVar_p^2 \1_{\{\tau_l \leq p\}}).
  \end{equation}
Summing for $p$ from $1$ to $l$, we obtain:
  \begin{align*}
  | \E( \Pi_k(t) \overline{\Pi_l(t)}) -1 | & \leq \frac{t^2 \exp(2Mt^2)
  }{B_k B_l}\ \E\left( \sum_{p=\tau_l}^l \RandVar_p^2 \right)
  \\&
  \leq \frac{t^2 \exp(2Mt^2)
  }{B_k B_l} \ \E\left( \RandVar_{\tau_l}^2 + \sum_{p=\tau_l+1}^l \RandVar_p^2 \right)
  \\ & \leq
   \frac{t^2 \exp(2Mt^2)
  }{B_k B_l}\ \E\left( \max_{1\leq j \leq l} \RandVar_j^2 + 2MB_l^2 \right).
  \end{align*}
Since $\E( \max_{1\leq j\leq l} \RandVar_j^2 /B_l^2)$ is uniformly bounded
by assumption, we obtain finally:
  \begin{equation}
  | \E( \Pi_k(t) \overline{\Pi_l(t)}) -1 | \leq C \frac{B_l}{B_k}\cdot
  \end{equation}
Write $\Pi_n^b(t)=\frac{1}{\sum_{k=1}^n b_k} \sum_{k=1}^n
b_k\Pi_k(t)$. Then
  \begin{align*}
  \E( |\Pi_n^b(t)-1|^2) &\leq \frac{2}{\left(\sum_{k=1}^n b_k\right)^2}
  \sum_{k=1}^n b_k \sum_{l=1}^k b_l | \E(\Pi_k(t) \overline{\Pi_l(t)})
  -1|\\
  & \leq  \frac{C}{\left(\sum_{k=1}^n b_k\right)^2}
  \sum_{k=1}^n b_k \sum_{l=1}^k b_l \frac{B_l}{B_k}\cdot
  \end{align*}
By assumption, $b_l B_l \leq C(B_l -B_{l-1})$. Summing from $1$ to
$k$, we get $\sum_{l=1}^k b_l B_l \leq C B_k$. Finally,
  \begin{equation}
  \label{Pi_n_b_sommable}
   \E( |\Pi_n^b(t)-1|^2)\leq \frac{C}{\left(\sum_{k=1}^n b_k\right)^2}
  \sum_{k=1}^n b_k \frac{B_k}{B_k} \leq \frac{C}{\sum_{k=1}^n b_k}\cdot
  \end{equation}
Since $b_k$ is bounded and $\sum b_k=+\infty$, there exists a
sequence $u_n$ such that $\sum_{k=1}^{u_n} b_k - n^2=O(1)$. Equation
\eqref{Pi_n_b_sommable} shows that $\E(| \Pi_{u_n}^b(t)-1|^2)$ is
summable. In particular, for almost every $\omega$, $\Pi_{u_n}^b(t)$
converges to $1$.

Consider now an arbitrary $m$, and choose $n$ with $u_n\leq
m<u_{n+1}$. Then
  \begin{equation}
  \Pi_{m}^b(t)=\frac{\sum_{k=1}^{u_n} b_k}{\sum_{k=1}^m b_k}
  \ \Pi_{u_n}^b(t)+\frac{1}{\sum_{k=1}^m b_k} \sum_{k=u_n+1}^m b_k
  \Pi_k(t).
  \end{equation}
Since $\frac{\sum_{k=1}^{u_n} b_k}{\sum_{k=1}^m b_k} \to 1$,
$\Pi_{u_n}^b(t) \to 1$, $\Pi_k(t)$ is bounded and
$\frac{\sum_{k=u_n+1}^m b_k}{\sum_{k=1}^m b_k} \to 0$, this shows
that $\Pi_m^b(t)$ converges to $1$. Hence, the third term of
\eqref{decomposition_en_trois} tends almost surely to $0$. This
concludes the proof.
\end{proof}

\subsection{Dynamical application}

\begin{proof}[\MakeThmProofTitle{\ref{TCLps_martingales}}]
Under the assumptions of Gordin's Theorem, there exist two functions
$g,h\in L^2$ such that $f=g-g\circ T+h$, and the sequence $h\circ
T^n$ is a reverse martingale difference for the filtration
$\boF_n=T^{-n}\boF$, i.e., $h\circ T^n$ is $\boF_n$-measurable and
$\E(h\circ T^n|\boF_{n+1})=0$. Moreover, a variance computation using
\eqref{eq:HypGordin} shows that $S_n f/\sqrt{n}$ is bounded in $L^2$.

Let us first show that $S_n f/\sqrt{n}$ has tight maxima. Since $S_n
h/\sqrt{n}= S_n f/\sqrt{n} + (g-g\circ T^n)/\sqrt{n}$, the sequence
$S_n h/\sqrt{n}$ is bounded in $L^2$, and therefore in $L^1$. Example
\ref{ex:TightMaxima} thus shows that $S_n h/\sqrt{n}$ has tight
maxima. To conclude the proof, we just have to check that $(g-g\circ
T^n)/\sqrt{n}$ has tight maxima. We have
  \begin{equation}
  \frac{\max_{1\leq k \leq n} | g-g\circ T^k|}{\sqrt{n}}
  \leq \frac{|g|}{\sqrt{n}} + \frac{ \max_{1\leq k\leq n}|g\circ
  T^k|}{\sqrt{n}}\cdot
  \end{equation}
Moreover, for any $c>0$,
  \begin{multline*}
  \P\left\{\max_{1\leq k\leq n}|g\circ
  T^k| /\sqrt{n} \geq c\right\} \leq \P\left\{ \max_{1\leq k \leq n} g^2\circ T^k /n
  \geq c^2\right\}\\ \leq c^{-2} \E( S_n g^2/n)=c^{-2} \E(g^2).
  \end{multline*}
Hence, this sequence is also tight. This concludes the proof of the
tightness of maxima of $S_n f/\sqrt{n}$.

Let us now turn to the proof of the almost sure central limit theorem.
Set $b_k=\frac{1}{k}$, $B_n=\sqrt{n}$, $\zeta=\int h^2\dd m$ and
$\RandVar_n=h\circ T^{n+1}$. We check the assumptions of Theorem
\ref{TCLps_martingale}. Birkhoff's ergodic Theorem applied to $h^2$ shows
that $h^2\circ T^n=o(n)$, hence $\RandVar_n/\sqrt{n}\to 0$ almost
everywhere. Moreover, $\sum_{k=1}^n \RandVar_j^2/n$ tends almost everywhere
to $\E(h^2)$, and the second condition of Theorem
\ref{TCLps_martingale} follows. The third condition is trivial.
Finally,
  \begin{equation}
  \E\bigl[ \max_{1\leq j\leq k} |\RandVar_j|^2/B_k^2 \bigr]
  \leq \E( S_k h^2 /k)= \E(h^2)
  \end{equation}
hence this sequence is bounded.

Therefore, Theorem~\ref{TCLps_martingale} applies and proves that,
almost everywhere,
  \begin{equation}
  \frac{1}{\log N} \sum_{k=1}^N \frac{1}{k} \gdelta_{S_k h(x)/
  \sqrt{k}} \lawto \boN(0, \E(h^2))\,.
  \end{equation}
Moreover, $S_k f=S_k h + g-g\circ T^k$. Again by Birkhoff's ergodic Theorem
applied to $g^2$, $(g-g\circ T^k)/\sqrt{k}$ tends almost surely to
$0$. The result follows.
\end{proof}

\bibliography{biblio}
\bibliographystyle{alpha}

\end{document}